%common header stuff

\documentstyle{amsppt}
\baselineskip18pt
\magnification=\magstep1
%\NoPageNumbers
%\NoRunningHeads
%\pagewidth{4.5in}
%\pageheight{7.0in}
\pagewidth{30pc}
\pageheight{45pc}
\hyphenation{co-deter-min-ant co-deter-min-ants pa-ra-met-rised
pre-print pro-pa-gat-ing pro-pa-gate
fel-low-ship Cox-et-er dis-trib-ut-ive}
\def\leaderfill{\leaders\hbox to 1em{\hss.\hss}\hfill}
\def\A{{\Cal A}}

\def\H{{\Cal H}}
\def\L{{\Cal L}}

\def\idest{i.e.,\ }
\def\wh{\widehat}
\def\ti{\widetilde}

\def\d{{\delta}}
\def\De{{\Delta}}

\def\th{{\theta}}
\def\i{{\iota}}
\def\k{{\kappa}}
\def\l{{\lambda}}

\def\te{\widetilde t}

\def\b0{\text{\bf 0}}

\def\ra{{\ \longrightarrow \ }}

\def\lan{{\langle}}
\def\ran{{\rangle}}

\def\dt{{\Bbb D \Bbb T}}

\def\zed{{\Bbb Z}}
\def\kyu{{\Bbb Q}}
\def\enn{{\Bbb N}}

\def\boxit#1{\vbox{\hrule\hbox{\vrule \kern3pt
\vbox{\kern3pt\hbox{#1}\kern3pt}\kern3pt\vrule}\hrule}}
\def\rabbit{\vbox{\hbox{\kern0pt
\vbox{\kern0pt{\hbox{---}}\kern3.5pt}}}}

\def\tableau#1{
        \hbox {
                \hskip -10pt plus0pt minus0pt
                \raise\baselineskip\hbox{
                \offinterlineskip
                \hbox{#1}}
                \hskip0.25em
        }
}

\def\tabCol#1{
\hbox{\vtop{\hrule
\halign{\strut\vrule\hskip0.5em##\hskip0.5em\hfill\vrule\cr\lower0pt
\hbox\bgroup$#1$\egroup \cr}
\hrule
} } \hskip -10.5pt plus0pt minus0pt}

\def\CR{
        $\egroup\cr
        \noalign{\hrule}
        \lower0pt\hbox\bgroup$
}

% Set up the map arrows for commutative diagrams.
%\def\mapright#1{\smash{
%     \mathop{\longrightarrow}\limits^{#1}}}

%Set up macro for commutative diagrams etc. (see Ex. 18.46 in TeXbook)

\def\blank#1#2{
%\hbox to {{#1}}{\vbox to {{#2}}}
\hbox to #1{\hfill \vbox to #2{\vfill}}
}

%Ross's table macros

\def\strut{\vrule height10pt depth5pt width0pt}

\topmatter
\title Decorated tangles and canonical bases
\endtitle

\author R.M. Green \endauthor
\affil 
Department of Mathematics and Statistics\\ Lancaster University\\
Lancaster LA1 4YF\\ England\\
{\it  E-mail:} r.m.green\@lancaster.ac.uk\\
\endaffil

\abstract
We study the combinatorics of fully commutative elements in Coxeter
groups of type $H_n$ for any $n > 2$.  Using the results, we construct
certain canonical bases for non-simply-laced generalized 
Tem\-per\-ley--Lieb algebras 
and show how to relate them to morphisms in the category of decorated 
tangles.
\endabstract

\thanks
The author was supported in part by an award from the Nuffield
Foundation.
\endthanks


\endtopmatter

\centerline{\bf To appear in the Journal of Algebra}

%Final version, 8th August 2001

\head Introduction \endhead

The category of decorated tangles was introduced by the author \cite{{\bf 4}}
using ideas of Martin and Saleur \cite{{\bf 12}}.  This allows a
generalization of the well-known diagram calculus \cite{{\bf 13}} for the 
Temperley--Lieb algebra to Coxeter systems of other types.  We showed
in \cite{{\bf 4}} and
\cite{{\bf 5}} how the endomorphisms in the category of decorated tangles
may be used to construct faithful representations of generalized
Temperley--Lieb algebras arising from Coxeter systems of types $B$,
$D$ or $H$.  By realizing these algebras, which are quotients of Hecke
algebras, in terms of the category, we can prove results
about their representation theory and structure which may otherwise
not be obvious, particularly in the case of type $H$.

In \cite{{\bf 6}}, it was shown that a generalized Temperley--Lieb algebra
arising from a Coxeter system of arbitrary type admits a canonical
basis (IC basis).  Formally, this is similar to the basis 
$\{C'_w : w \in W\}$ for the Hecke algebra introduced by Kazhdan and 
Lusztig \cite{{\bf 10}}, although the precise relationship between the two 
is not completely obvious.

The goal of this paper is to tie these two theories together by
showing how decorated tangles may be used to
describe certain canonical bases for generalized Temperley--Lieb
algebras.  This is easily done in types $D$ (and $A$) by using \cite{{\bf 6},
Theorem 3.6} and \cite{{\bf 4}, Theorem 4.2}, so we concentrate here on
the nontrivial cases of Coxeter systems of types $B$ and $H$.

According to \cite{{\bf 1}}, interesting algebras and representations
defined over $\enn$ come from category theory, and are best
understood when their categorical origin has been discovered.  It is
conjectured \cite{{\bf 7}, Conjecture 1.2.4} that the canonical basis in
the preceding paragraph should have structure constants in $\enn[v,
v^{-1}]$, and the results of this paper may be regarded as discovering
the categorical origin of this phenomenon in certain special cases.

Further motivation for our results is as follows.  In type $B$, it is known
from \cite{{\bf 7}, Theorem 2.2.1} that the canonical basis is the 
projection of a certain subset of the Kazhdan--Lusztig basis, but no
purely combinatorial construction is given; we give such a
construction here (Theorem 2.2.5).
Conversely, in type $H$, we have a combinatorial construction of a basis
with very convenient and useful properties \cite{{\bf 5}, \S4}, but we
have no simple abstract characterisation of the basis; we show here that the
basis of decorated tangles from \cite{{\bf 5}} is identical to the
canonical basis (Theorem 2.1.3).

Our methods of proof are combinatorial, and a by-product of the proof
is a construction of the canonical basis in types $B$ and $H$ which
avoids diagrams and Kazhdan--Lusztig theory and uses only the 
combinatorics of fully commutative words (see Theorem 3.4.3 and
Theorem 5.2.1).  It is therefore necessary to develop an understanding
of this combinatorics in type $H$, which we do in \S3.  As we see in
\S5, the combinatorics in type $H$ behaves
rather like a superset of the combinatorics of fully commutative
expressions in type $B$; the latter has been studied extensively by
Stembridge \cite{{\bf 15}, \S\S2--6} from a very different perspective.

Apart from proving the claim in \cite{{\bf 6}, Remark 2.4 (1)}, the
results of this paper are of interest because of their applications.
For example, our results here are related in a very precise way to 
Jones' planar algebras \cite{{\bf 9}}; this is the subject of a future paper.
Furthermore, there are applications to Kazhdan--Lusztig theory: the
curious combinatorial rules in definitions 2.1.1 and 2.2.4 can be
interpreted to give a combinatorially explicit description of those
Kazhdan--Lusztig cells that contain fully commutative elements, certainly
when the associated Coxeter groups are finite, and conjecturally in general.

\head 1. Preliminaries \endhead

\subhead 1.1 Decorated tangles \endsubhead

We start by recalling the definition of the category of decorated
tangles which was introduced in \cite{{\bf 4}}.  This is a mild
generalization of Martin and Saleur's diagram calculus for the blob
algebra in \cite{{\bf 12}}.

\definition{Definition 1.1.1}
A tangle is a portion of a knot diagram contained in a rectangle.  The
tangle is incident with the boundary of the rectangle only on the
north and south faces, where it intersects transversely.  The
intersections in the north (respectively, south) face are numbered
consecutively starting with node number $1$ at the western (\idest the
leftmost) end.

Two tangles are {\it equivalent} if there exists an isotopy of the
plane carrying
one to the other such that the corresponding faces of the rectangle
are preserved setwise.  Two tangles are {\it vertically equivalent} if they
are equivalent in the above sense by an isotopy which preserves setwise each
vertical cross-section of the rectangle.
\enddefinition

We call the edges of the rectangular frame ``faces'' to avoid
confusion with the ``edges'' which are the arcs of the tangle.

\definition{Definition 1.1.2}
A decorated tangle is a crossing-free tangle in which an arc exposed
to the west face of the rectangular frame is allowed to carry decorations.
(Any arc not exposed to the west face 
of the rectangular frame must be undecorated.)
\enddefinition

By default, the decorations will be discs, although we will need two
kinds of decorations for some of our later applications.

Any decorated tangle consists only of loops and edges,
none of which intersect each other.  A typical example is shown in
Figure 1.

\topcaption{Figure 1} A decorated tangle \endcaption
\centerline{
\hbox to 3.027in{
\vbox to 0.888in{\vfill
        \includegraphics{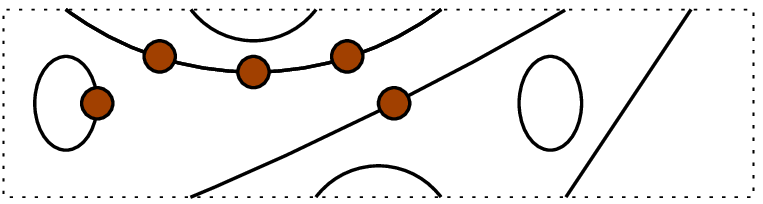}
}
\hfill}
}

\definition{Definition 1.1.3}
The category of decorated tangles, $\dt$, has as its objects the
natural numbers.  The morphisms from $n$ to $m$
are the equivalence classes of decorated tangles with $n$ nodes in the 
north face and $m$ in the south.  The
source of a morphism is the number of points in the north face of the
bounding rectangle, and the target is the number of points in the
south face.  Composition of morphisms works by concatenation of the
tangles, matching the relevant south and north faces together.
\enddefinition

Note that for there to be any morphisms from $n$ to $m$, it is
necessary that $n+m$ be even.

We now define the algebra of decorated tangles.

\definition{Definition 1.1.4}
Let $R$ be a commutative ring with $1$ and let 
$n$ be a positive integer.  Then the $R$-algebra $\dt_n$ has as a
free $R$-basis the morphisms from $n$ to $n$, where the multiplication
is given by the composition in $\dt$.
\enddefinition

\definition{Definition 1.1.5}
The edges in a tangle $T$ which connect nodes (\idest not the loops)
may be classified
into two kinds: propagating edges, which link a node in the north
face with a node in the south face, and non-propagating edges, which
link two nodes in the north face or two nodes in the south face.
\enddefinition

\subhead 1.2 Canonical bases \endsubhead

Let $X$ be a Coxeter graph, of arbitrary type,
and let $W(X)$ be the associated Coxeter group with distinguished
set of generating involutions $S(X)$.  Denote by $\H(X)$ the Hecke
algebra associated to $W(X)$.  (A good reference for Hecke algebras is
\cite{{\bf 8}, \S7}.)  Let $\A$ be the ring of
Laurent polynomials, $\zed[v, v^{-1}]$.  The $\A$-algebra $\H(X)$ has 
a basis consisting of elements $T_w$, with $w$ ranging over $W(X)$, 
that satisfy $$T_s T_w = 
\cases
T_{sw} & \text{ if } \ell(sw) > \ell(w),\cr
q T_{sw} + (q-1) T_w & \text{ if } \ell(sw) < \ell(w),\cr
\endcases$$ where $\ell$ is the length function on the Coxeter group
$W(X)$, $w \in W(X)$, and $s \in S(X)$.
The parameter $q$ is equal to $v^2$.

Let $J(X)$ be the two-sided ideal of $\H(X)$ generated by the elements $$
\sum_{w \in \lan s, s' \ran} T_w,
$$ where $(s, s')$ runs over all pairs of elements of $S(X)$
that correspond to adjacent nodes in the Coxeter graph.  
(If the nodes corresponding to $(s, s')$ are connected by a
bond of infinite strength, then we omit the corresponding relation.)

\definition{Definition 1.2.1}
Following Graham \cite{{\bf 3}, Definition 6.1}, we define the generalized
Temperley--Lieb algebra $TL(X)$ to be
the quotient $\A$-algebra $\H(X)/J(X)$.  We denote the corresponding
epimorphism of algebras by $\th : \H(X) \ra TL(X)$.
\enddefinition

\definition{Definition 1.2.2}
A product $w_1w_2\cdots w_n$ of elements $w_i\in W(X)$ is called
{\it reduced} if $\ell(w_1w_2\cdots w_n)=\sum_i\ell(w_i)$.  We reserve
the terminology {\it reduced expression} for reduced products 
$w_1w_2\cdots w_n$ in which every $w_i \in S(X)$.

Call an element $w \in W(X)$ {\it complex} if it can be written 
as a reduced product $x_1 w_{ss'} x_2$, where $x_1, x_2 \in W(X)$ and
$w_{ss'}$ is the longest element of some rank 2 parabolic subgroup 
$\lan s, s'\ran$ such that $s$ and $s'$ correspond to adjacent nodes
in the Coxeter graph.

Denote by $W_c(X)$ the set of all elements of $W(X)$
that are not complex.  The {\it fully commutative} elements of $W(X)$ are
defined in \cite{{\bf 14}} to be those elements $w \in W$ for which any reduced
expression can be transformed into any other by repeatedly applying
short braid relations, \idest by iterated commutation
of pairs of consecutive generators in the expression.  The elements
$W_c(X)$ are precisely the fully commutative elements of $W(X)$ by
\cite{{\bf 14}, Proposition 1.1}.

We define the {\it content} of $w\in W$ to be the set $c(w)$
of Coxeter generators $s\in S$ that appear in a reduced expression for 
$w$.  This can be shown not to depend on the reduced expression chosen
by using the theory of Coxeter groups, e.g. by applying 
Matsumoto's Theorem.

Let $t_w$ denote the image of the basis element $T_w \in \H(X)$ in
the quotient $TL(X)$.
\enddefinition

\proclaim{Theorem 1.2.3 (Graham)}
The set $\{ t_w : w \in W_c \}$ 
is an $\A$-basis for the algebra $TL(X)$.  \endproclaim

\demo{Proof}
See \cite{{\bf 3}, Theorem 6.2}.
\qed\enddemo

\definition{Definition 1.2.4}
For each $s\in S(X)$, we define $b_s=v^{-1}t_s+v^{-1}t_e$, where $e$
is the identity element of $W$.  If 
$w\in W_c(X)$ and $s_1s_2\cdots s_n$ is a reduced expression
for $w$, then we define $b_w = b_{s_1}b_{s_2}\cdots b_{s_n}$.  We take 
the empty product $b_e$ to be the identity element $t_e$ of $TL(X)$.

Note that $b_w$ does not depend on the choice of reduced
expression for $w$.

It is known that $\{b_w : w\in W_c\}$ is a basis for the $\A$-module
$TL(X)$; this can be deduced from Theorem 1.2.3 by using \cite{{\bf 6},
Lemma 1.5}.  We shall call it the {\it monomial basis}.
\enddefinition

We now recall a principal result of \cite{{\bf 6}}, which establishes
the canonical basis for $TL(X)$.  This basis is a direct analogue of 
the important Kazhdan--Lusztig basis of the Hecke algebra $\H(X)$.
  
Fix a Coxeter graph, $X$.  Let $I=W_c(X)$, 
let $\A^- = \zed[v^{-1}]$, and let $\,\bar{\ }\,$ be the
involution on the ring $\A = {\Bbb Z}[v, v^{-1}]$ which satisfies 
$\bar{v} = v^{-1}$.

By \cite{{\bf 6}, Lemma 1.4}, the algebra $TL(X)$ has a $\zed$-linear 
automorphism of order $2$ that sends $v$ to $v^{-1}$ and $t_w$ to 
$t_{w^{-1}}^{-1}$.  We denote this map also by $\,\bar{\ }\,$.

Let $\L$ be the free $\A^-$-submodule of $TL(X)$ with basis
$\{\te_w : w \in W_c\}$, where $\te_w := v^{-\ell(w)} t_w$,
and let $\pi : \L \ra \L/v^{-1}\L$ be the canonical projection.

\proclaim{Theorem 1.2.5}
There exists a unique basis $\{ c_w : w \in W_c\}$ for $\L$
such that $\overline{c_w} = c_w$ and $\pi(c_w) = \pi(\te_w)$
for all $w\in W_c$.
\endproclaim

\demo{Proof}
This is \cite{{\bf 6}, Theorem 2.3}.
\qed\enddemo

The basis $\{c_w : w \in W_c\}$ is called the {\it IC basis} (or the 
{\it canonical basis}) of $TL(X)$.  It depends on the 
$t$-basis, the involution $\,\bar{\ }\,$, and the lattice $\L$.  

\head 2. Main results \endhead

\subhead 2.1 Type $H$ \endsubhead

A Coxeter system of type $H_n$ is given by the Coxeter graph in Figure
2.

\topcaption{Figure 2} Coxeter graph of type $H_n$ \endcaption
\centerline{
\hbox to 3.138in{
\vbox to 0.402in{\vfill
        \includegraphics{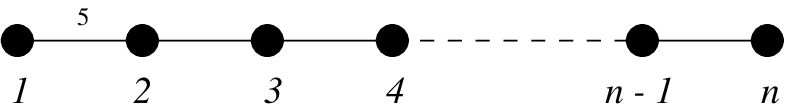}
}
\hfill}
}

The algebras $TL(H_n)$ were first studied by Graham in \cite{{\bf 3}, \S7}.
The basic properties of these algebras are also described in
\cite{{\bf 5}}, to which the reader is referred for further explanatory
comments and examples.  We summarise some of the more important
properties here.

The unital algebras $TL(H_n)$ 
over $\A$ can be defined by generators $b_1, b_2, \ldots b_n$ (as in
Definition 1.2.4) and relations $$\eqalign{
b_i^2 &= \d b_i, \cr
b_i b_j &= b_j b_i \text{\quad if \ $|i - j| > 1$},\cr
b_i b_j b_i &= b_i \text{\quad if \ $|i - j| = 1$ \ and \ $i, j > 1$},\cr
b_i b_j b_i b_j b_i &= 3 b_i b_j b_i - b_i
\text{\quad if \ $\{i, j\} = \{1, 2\}$},\cr
}$$ where $\d := [2] = v + v^{-1}$.
All the algebras $TL(H_n)$ are finite dimensional, even though the
Coxeter group $W(H_n)$ is infinite for $n > 4$.  

The algebras $TL(H_n)$ can also be described in terms of decorated
tangles; for this, we need the concept of an $H$-admissible diagram
from \cite{{\bf 5}, Definition 2.2.1}.

\definition{Definition 2.1.1}
An $H$-admissible diagram with $n$ edges is an element of $\dt_n$
with no loops which satisfies the following conditions.

\item{\rm (i)}
{No edge may be decorated if all the edges are propagating.}
\item{\rm (ii)}
{If there are non-propagating edges in the diagram, then either there
is a decorated edge in the north face connecting nodes 1 and 2, or
there is a non-decorated edge in the north face connecting nodes $i$
and $i+1$ for some $i > 1$.
A similar condition holds for the south face.}
\item{\rm(iii)}
{Each edge carries at most one decoration.}

The $H$-admissible diagram $U_i$, where $1 \leq i \leq n$, is the diagram
all of whose edges are propagating and undecorated, except for those
attached to nodes $i$ and $i+1$ in the north face, and nodes $i$ and
$i+1$ in the south face.  These four nodes are connected in the pairs
given, using decorated edges if $i = 1$, and using undecorated edges
if $i > 1$.
\enddefinition

The main result of \cite{{\bf 5}} is the following.

\proclaim{Theorem 2.1.2}
Let $\De_n^H$ be the $\A$-algebra with
basis given by the $H$-admissible diagrams with $n$ edges and
multiplication induced from the multiplication on $\dt_n$, subject to
the reduction rules in Figure 3.  Then the map sending $b_i$ to $U_i$
extends to an isomorphism of $\A$-algebras between $TL(H_{n-1})$
and $\De_n^H$.
\endproclaim

\vfill\eject

\topcaption{Figure 3} Reduction rules for $\Delta_n$ \endcaption
\centerline{
\hbox to 1.041in{
\vbox to 1.916in{\vfill
        \includegraphics{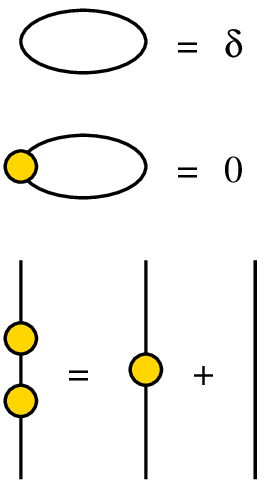}
}
\hfill}
}

\demo{Proof}
This is \cite{{\bf 5}, Theorem 3.4.2}.
\qed\enddemo

Recall \cite{{\bf 5}, Remark 2.2.3} 
that the first two relations in Figure 3 determine which scalar
to multiply by when a loop is removed, and the third relation
expresses a diagram as a sum of two other diagrams with fewer decorations.

One of the two main results of this paper is the next result.

\proclaim{Theorem 2.1.3}
The diagram basis described in Theorem 2.1.2 is the
canonical basis of $TL(H_{n-1})$ in the sense of Theorem 1.2.5.
\endproclaim

\subhead 2.2 Type B \endsubhead

We have results analogous to those in \S2.1 for Coxeter systems of
type $B$.  Here, nodes 1 and 2 in the Coxeter graph are connected by a
bond of strength 4 rather than 5.  

The finite dimensional, unital algebras $TL(B_n)$ 
over $\A$ can be defined by generators $b_1, b_2, \ldots, b_n$
and relations $$\eqalign{
b_i^2 &= \d b_i, \cr
b_i b_j &= b_j b_i \text{\quad if \ $|i - j| > 1$},\cr
b_i b_j b_i &= b_i \text{\quad if \ $|i - j| = 1$ \ and \ $i, j > 1$},\cr
b_i b_j b_i b_j  &= 2 b_i b_j
\text{\quad if \ $\{i, j\} = \{1, 2\}$},\cr
}$$ where $\d := [2] = v + v^{-1}$.

We now describe $TL(B_n)$ in terms of decorated tangles.  We first
recall the description from \cite{{\bf 4}}, and then show how it can be
adapted to give a result analogous to Theorem 2.1.3.  It is convenient
for this purpose to introduce some further terminology to classify tangles.

\definition{Definition 2.2.1}
We say an edge in a tangle is of type $p_1$ if it connects node $1$ in
the north face to node $1$ in the south face.  We say the edge is of
type $p_3$ if it does not involve either node $1$ in the north face or
node $1$ in the south face.  Otherwise, we say the edge is of type
$p_2$.
\enddefinition

%2 and 3 were originally the wrong way round above

The following definition is compatible with that in \cite{{\bf 4}, Theorem 4.1}.

\definition{Definition 2.2.2}
A $B$-admissible diagram with $n$ edges is an element of $\dt_n$ with
no loops which satisfies one of the following three mutually
exclusive conditions:

\item{\rm (B$1$)}
{There is an undecorated edge of type $p_1$ (and no other decorations).}
\item{\rm (B$1'$)}
{There is a decorated edge of type $p_1$ 
and at least one non-propagating edge.}
\item{\rm (B$2$)}
{There are two decorated edges of type $p_2$.}
\enddefinition

We recall one of the main results of \cite{{\bf 4}}, noting that the
diagrams $U_i$ are $B$-admissible.

\proclaim{Theorem 2.2.3}
Let $\De_n^B$ be the $\kyu[v, v^{-1}]$-algebra with
basis given by the $B$-ad\-miss\-ible diagrams with $n$ edges and
multiplication induced from the multiplication on $\dt_n$, subject to
the reduction rules in Figure 4.  Then the map sending $b_1$ to $2U_1$
and $b_i$ to $U_i$, for $i > 1$,
extends to an isomorphism of $\kyu[v, v^{-1}]$-algebras between $TL(B_{n-1})$
and $\De_n^B$.
\endproclaim

\vfill\eject

\topcaption{Figure 4} Reduction rules for type $B$ \endcaption
\centerline{
\hbox to 1.222in{
\vbox to 1.916in{\vfill
        \includegraphics{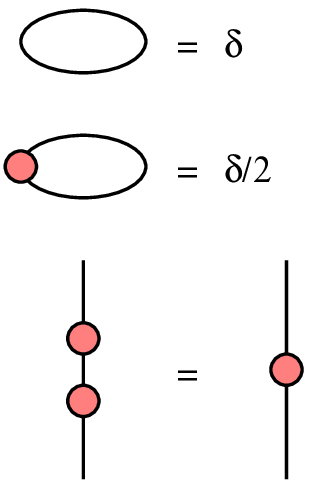}
}
\hfill}
}

\demo{Proof}
See \cite{{\bf 4}, Theorem 4.1}.
\qed\enddemo

The reason for the presence of $\kyu$ rather than $\zed$ in Theorem 2.2.3 is
because our integral form of the algebra $TL(B_{n-1})$ is not the same
as the integral form used in \cite{{\bf 4}}.  We will return to this issue
in \S5.1.

The basis in Theorem 2.2.3 is not the canonical basis, although it is
closely related to it.  In order to
describe the canonical basis combinatorially, it is convenient to introduce a
second type of decoration, denoted by a square.  This satisfies the
relations in Figure 5.  The first rule given is the definition of
the square decoration, and the others are reduction rules which are
immediate consequences of the definition and the relations in Figure
4; we shall use these freely without further comment.

\vfill\eject

\topcaption{Figure 5} The square decoration \endcaption
\centerline{
\hbox to 1.916in{
\vbox to 2.430in{\vfill
        \includegraphics{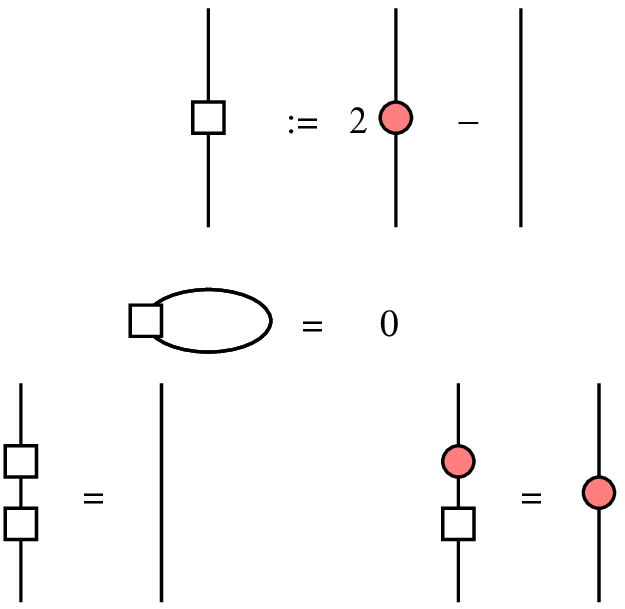}
}
\hfill}
}

\definition{Definition 2.2.4}
A $B$-canonical diagram with $n$ edges is an element of $\dt_n$ with
no loops that satisfies one of the following three mutually
exclusive conditions C$1$, C$1'$, C$2$ below.  All the decorations are
square, except for those on the two edges involved in C$2$.

\item{\rm (C$1$)}
{There is an undecorated edge of type $p_1$ (and no other decorations).}
\item{\rm (C$1'$)}
{There is an edge of type $p_1$ with a square decoration,
and at least one non-propagating edge.}
\item{\rm (C$2$)}
{There are two edges of type $p_2$ with circular decorations.}
\enddefinition

We denote by $C_n^1$ the set of all elements of $\dt_n$ that
satisfy C$1$ or C$1'$, and we denote by $C_n^2$ the set of elements
of $\dt_n$ that satisfy C$2$.  The example in Figure 6 is an element
of $C_8^2$.

\topcaption{Figure 6} A $B$-canonical diagram \endcaption
\centerline{
\hbox to 3.638in{
\vbox to 0.888in{\vfill
        \includegraphics{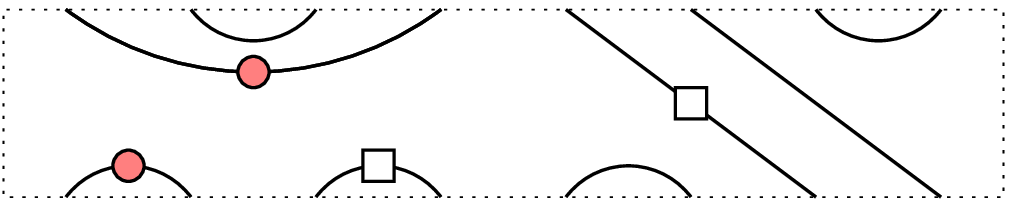}
}
\hfill}
}

Our second main result is the following.

\proclaim{Theorem 2.2.5}
The set $C_n := C_n^1 \cup \{2D: D \in C_n^2\}$ is the canonical basis 
of $TL(B_{n-1})$ in the sense of Theorem 1.2.5.
\endproclaim

\subhead 3. Combinatorial construction of the basis in type $H$ \endsubhead

The aim of \S3 is to obtain an explicit, elementary description of the 
canonical basis in type $H$ in terms of the monomial basis $\{b_w : w
\in W_c\}$ (see Definition 1.2.4).
The answer turns out to be more complicated than in the finite
dimensional simply-laced case, where the basis may be described in
terms of certain monomials \cite{{\bf 6}, Theorem 3.6}.

In \S3, all computations take place in $TL(H_{n-1})$ over the ring
$\A$ unless otherwise stated.

\subhead 3.1 Deletion properties for monomials in type $H$ \endsubhead

\proclaim{Lemma 3.1.1}
The $\zed$-span of the monomial basis of $TL(H_{n-1})$ is equal to
the $\zed$-span of the diagram basis of $TL(H_{n-1})$ described in
Theorem 2.1.2.
\endproclaim

\demo{Proof}
For the purposes of the proof, we shift base rings and
regard $TL(H_{n-1})$ as the $\zed[\d]$-algebra generated by the 
monomial basis elements $b_i$.  (The generators and relations in \S2.1
and the relations in the diagram algebra show that this makes sense.)
Denote by $TL_{\d} = TL_{\d}(H_{n-1})$ the $\zed[\d]$-version of the algebra.

From Definition 1.2.4, we know that $\{b_w: w \in W_c\}$ is an
$\A$-basis for $TL(H_{n-1})$.  We claim that it is a
$\zed[\d]$-basis for the version of the algebra we consider here.  To
prove this, it suffices to show that the set given is a spanning set.
Let $a \in TL_\d$.  Then, since $a$ is fixed by $\bar{ }$ (since $\d$ and
$b_i$ are), we know that when $a$ is expressed as a linear combination
of $b_w$ (with coefficients in $\A$), the coefficients are fixed by
$\bar{ }$.  An easy induction on the degree of Laurent polynomials
shows that the elements of $\A$ fixed by
$\bar{ }$ are precisely the elements of $\zed[\d]$, which completes
the claim.

A similar argument shows that 
the diagram $\A$-basis elements lie in $TL_{\d}$, because they are
$\zed$-linear combinations of monomials in the $b_i$ by \cite{{\bf 5},
Proposition 3.2.8} and therefore fixed by $\bar{ }$.

Let $b_w = \sum_{x} \l_x a_x$, where the $a_x$ are diagram
basis elements.
To complete the proof, it suffices to show that all the $\l_x$ are
integers.  The reduction rules in Figure 3 show that this will
happen if no loops appear in the diagram corresponding to the monomial
$b_w$.  If $k$ loops appear, the rules show that $\l_x$ will equal
$\d^k z$ for some (possibly zero) integer $z$.  Since the $b_w$ form a
$\zed[\d]$-basis for $TL_{\d}$ 
and $\d$ is not invertible in $\zed[\d]$, we must have $k = 0$
in each case.  This completes the claim.
\qed\enddemo

The next definition will turn out (once we have proved Theorem 3.4.3) 
to agree with the lattice $\L$ of \S1.2.

\definition{Definition 3.1.2}
The lattice $\L_H$ is the free $\A^-$-module on the monomial basis (or
diagram basis) of $TL(H_{n-1})$.  Let $\pi_H : \L_H \ra
\L_H/v^{-1}\L_H$ be the canonical projection.
\enddefinition

The equivalence of the two versions of the definition is immediate
from Lemma 3.1.1.

Recall from Definition 1.2.4 that $b_e = \te_e$ and $b_i=\te_i+v^{-1}\te_e$.
The following result, which is analogous to \cite{{\bf 6}, Lemma 3.3}, 
is of key importance in the sequel.

\proclaim{Lemma 3.1.3}
Let $b = b_{i_1} b_{i_2} \cdots b_{i_k}$ be an arbitrary monomial in the
generators $b_i (1 \leq i < n)$.  Let $\wh{b}(l)$ denote the
monomial obtained from $b$ by omission of the $l$-th term, $b_{i_l}$.
Then if $b \in v^m \L_H$ for some integer $m$, we have $\wh{b}(l) \in
v^{m+1} \L_H$.
\endproclaim

\demo{Proof}
Consider the (not necessarily $H$-admissible) diagram corresponding to
the monomial $b$.  This has a certain number, $c$, of loops.  It is
clear from the nature of the diagrams that removing one generator
$b_{i_l}$ will change the number of loops by at most $1$, so the
number of loops, $\wh{c}$, in $\wh{b}(l)$ satisfies 
$\wh{c} \leq c
+ 1$.  Expanding $\wh{b}(l)$ in terms of the diagram basis and
applying Lemma 3.1.1 shows that $\wh{b}(l) \in v^{m+1} \L_H$, as required.
\qed\enddemo

It is convenient to introduce some terminology to describe the
individual entries in a monomial in the generators $b_i$.  We are
particularly interested in the case of monomials which correspond to
fully commutative reduced words (as in Definition 1.2.4), but we will
also need the terminology in more general situations.  We refer below 
to the instances of the individual generators $b_i$ in a monomial as 
``letters'', for the sake of clarity.

\definition{Definition 3.1.4}
Let $b = b_{i_1} b_{i_2} \cdots b_{i_k}$ be a monomial in the
generators $b_i$.  Let $p = i_l$ for some $l$.

We say the letter $b_p$ is {\it internal} if, after applying
commutations, the monomial $b$ may be transformed to
a monomial in which the (same occurrence of the) letter $b_p$ occurs 
as the middle term of a 
sequence $b_q b_p b_q$ for some letter $b_q$ that does not commute
with $b_p$.  Otherwise, we say the letter $b_p$ is {\it external}.

We say the letter $b_p$ is {\it lateral} (with respect to an internal
letter $b_q$) if it is external and, after
commutations, it can be moved adjacent to the internal letter $b_q$ to
form a subsequence $b_p b_q b_p$, where $b_p$ and $b_q$ do not commute.  We
say the letter $b_p$ is {\it bilateral} if it is lateral with respect
to two different internal letters.

We will also apply these definitions in the obvious way to reduced
expressions for fully commutative elements.
\enddefinition

\remark{Remark 3.1.5}
An important fact is that if $b = b_w$ is an element of the monomial
basis, then the only possible internal letters are occurrences of
$b_1$ and $b_2$.
\endremark

To illustrate these concepts, here are some examples.

\example{Example 3.1.6}
Consider the fully commutative reduced expression $$
b = b_{i_1} \cdots b_{i_7} := b_1 b_2 b_3 b_1 b_2 b_1 b_2
$$ in $TL(H_3)$.
The letters $b_{i_2}$, $b_{i_5}$ and $b_{i_6}$ are internal, and the
others are external.  Of the external letters, $b_{i_1}$, $b_{i_4}$
and $b_{i_7}$ are lateral, and of these, $b_{i_4}$ is the only
bilateral one.  
\endexample

It is important to identify certain configurations that are not
allowed to occur in the diagrammatic representation of a monomial
basis element.  The short dashed lines in
the figures indicate that the two curves are part of the same
generator $b_i$ in a tangle for the monomial.

\topcaption{Figure 7} Some impermissible configurations \endcaption
\centerline{
\hbox to 5.791in{
\vbox to 3.027in{\vfill
        \includegraphics{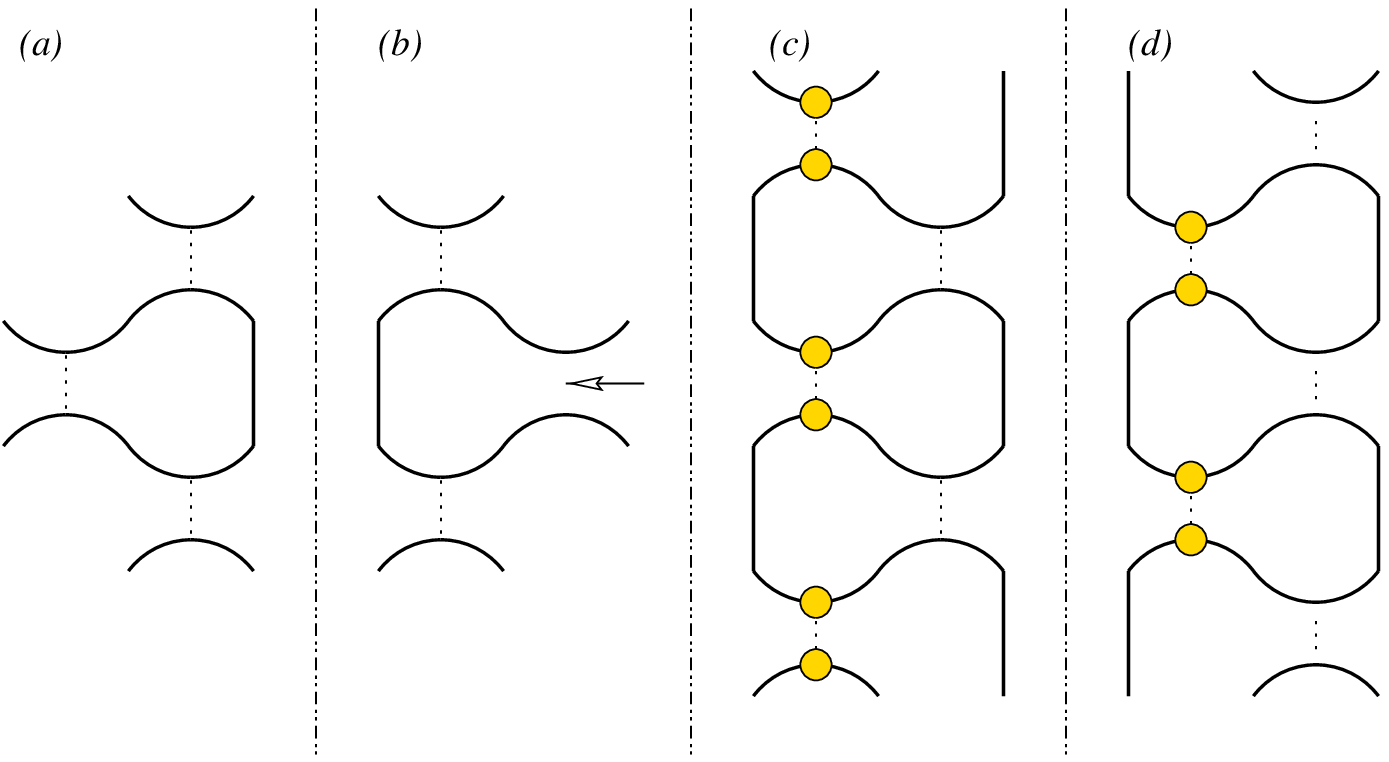}
}
\hfill}
}

\topcaption{Figure 8} Impermissible configuration (e) \endcaption
\centerline{
\hbox to 3.541in{
\vbox to 1.888in{\vfill
        \includegraphics{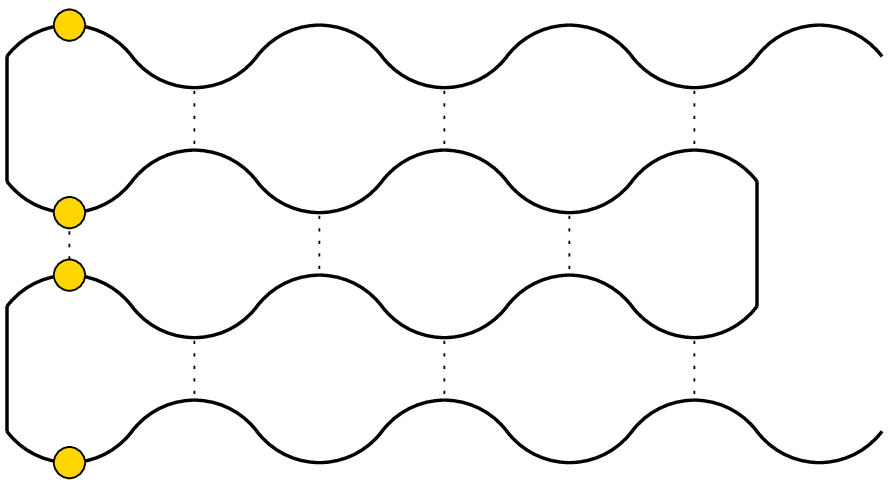}
}
\hfill}
}

\topcaption{Figure 9} Impermissible configuration (f) \endcaption
\centerline{
\hbox to 1.541in{
\vbox to 0.777in{\vfill
        \includegraphics{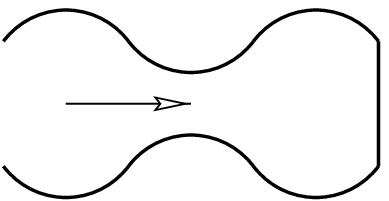}
}
\hfill}
}

\proclaim{Lemma 3.1.7}
Consider a tangle $T_V$ which is vertically equivalent to a monomial
$b_w$ for $w \in W_c$.  Then $T_V$ does not contain any segments
vertically equivalent to those in parts (a) to (d) of Figure 7.
\endproclaim

\demo{Proof}
Note first that commutation of
generators $b_i b_j = b_j b_i$ preserves vertical equivalence.
Parts (a), (c) and (d) of Figure 7 correspond, after applying
commutations if necessary, to sequences $b_{i+1} b_i b_{i+1}$ (for
some $i > 1$), $b_1b_2b_1b_2b_1$ and $b_2b_1b_2b_1b_2$, none of which
can occur in a fully commutative monomial by definition.  (Compare
with ``Property R'' depicted in Figure 2 of \cite{{\bf 2}}.)

For (b), either the configuration corresponds to a sequence $b_i
b_{i+1} b_i$ (again for $i > 1$) or there is a smaller configuration
of the same shape as (b) inserted in the position shown by the arrow.
Iterating this argument rightwards shows that there is an occurrence
of $b_i b_{i+1} b_i$ (for some $i > 1$) somewhere, which is not allowed.
\qed\enddemo

\remark{Remark 3.1.8}
Figure 8 shows an impossible situation which violates condition
(a) of Figure 7.  Similar counterexamples exist where the lobes of the
arc are either longer or shorter than the one shown in the diagram.

The situation in Figure 9 is also impermissible (compare condition (b)
of Figure 7).  In this situation, there is either a smaller configuration of
the same shape as (f) inserted in the position shown by the arrow, or
we have a violation of condition (a), (b) or (c) of Figure 7.

These facts are very important in the proof of the next result,
which is a more specialized version of Lemma 3.1.3.
\endremark

\proclaim{Proposition 3.1.9}
Let $b = b_w \in \L_H$ be an element of the monomial basis.  Let $b_p$ be a
letter in the monomial $b$, and let $\wh{b}$ be the monomial
obtained from $b$ by deleting the letter $b_p$.  Suppose $\wh{b} \not\in
\L_H$.  Then $w$ can be parsed in one of the following ways, where the
deleted letter is the barred one:
\item{\rm(i)}{$w = w_1 s_1 s_3 \cdots s_{2k-1} s_2 s_4 \cdots
\overline{s_{2k}} s_1 s_3 \cdots s_{2k-1} w_2$, where $p = 2k > 2$,}
\item{\rm(ii)}{$w = w_1 s_1 s_3 \cdots s_{2k+1} s_2 s_4 \cdots s_{2k}
s_1 s_3 \cdots s_{2k-1} w_2 \overline{s_{2k}} s_{2k+1} w_3$, where $p
= 2k$ and $s_{2k+1}$ commutes with every member of $c(w_2)$,}
\item{\rm(iii)}{$w = w_1 s_{2k+1} \overline{s_{2k}} w_2
s_1 s_3 \cdots s_{2k-1} s_2 s_4 \cdots s_{2k}
s_1 s_3 \cdots s_{2k+1} w_3$, where $p 
= 2k$ and $s_{2k+1}$ commutes with every member of $c(w_2)$,}
\item{\rm(iv)}{$w = w_1 s_q \overline{s_p} s_q w_2$, where 
$\{p, q\} = \{1, 2\}$.}

Conversely, if $w$ can be parsed in one of the ways above, removal of
the barred letter results in a monomial $\wh{b} \not\in \L_H$.
\endproclaim

\definition{Definition 3.1.10}
We say the letter $\overline{s_p}$ in $w$ is critical of type (i), type (ii)
or type (iii) respectively if,
after commutations, $w$ can be parsed as in parts (i), (ii) or (iii)
respectively of Proposition 3.1.9.  (If $\overline{s_p}$ satisfies condition
(iv), this is the same as saying $\overline{s_p}$ is internal.)
\enddefinition

\demo{Proof of Proposition 3.1.9}
We proceed using the diagram calculus, by considering an element $T_V$
vertically equivalent to the tangle corresponding to the monomial $b_w$.  

Recall that none of the configurations in figures 7 and 8 is allowed
to occur.  If the letter to be deleted is internal as in case (iv),
the situation is as shown in Figure 10.

\topcaption{Figure 10} Deleting a letter \endcaption
\centerline{
\hbox to 0.597in{
\vbox to 1.277in{\vfill
        \includegraphics{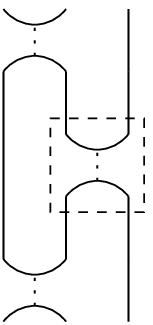}
}
\hfill}
}

Suppose we are not in case (iv).  
A routine, but nontrivial, case by case check using Lemma 3.1.7 and
Remark 3.1.8 shows that if the letter
being deleted forms a loop to the left, we must be in the situation shown in
Figure 11.  If, on the other hand, the letter being deleted forms a
loop to the right, we must be in the situation shown in Figure 12 or
its top-bottom mirror image.  (As in Figure 8, the lobes of the arcs
occurring may be longer or shorter than those shown.)

\topcaption{Figure 11} Deleting a critical letter of type (i) \endcaption
\centerline{
\hbox to 3.027in{
\vbox to 1.527in{\vfill
        \includegraphics{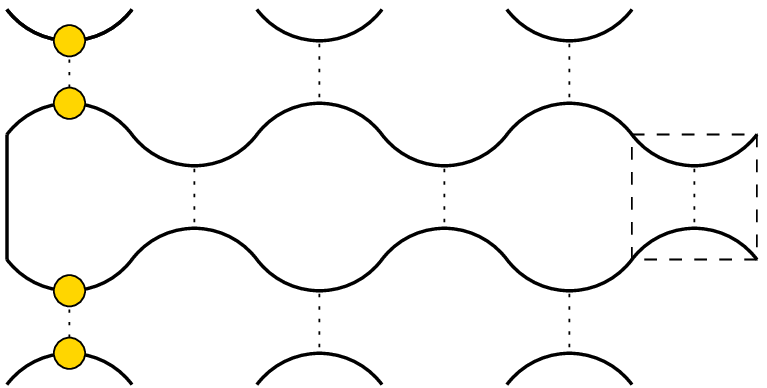}
}
\hfill}
}

\topcaption{Figure 12} Deleting a critical letter of type (ii) \endcaption
\centerline{
\hbox to 3.555in{
\vbox to 2.291in{\vfill
        \includegraphics{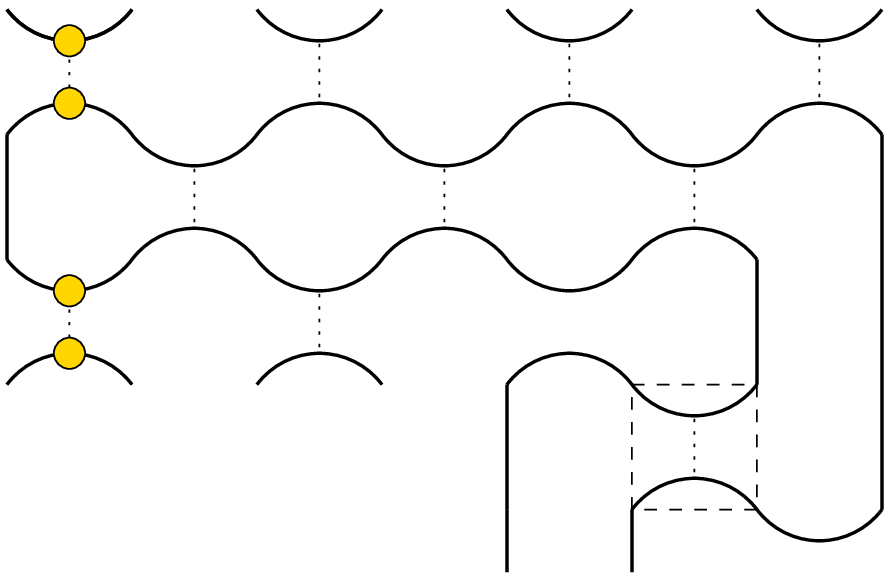}
}
\hfill}
}

Figure 11 corresponds to case (i), Figure 12 corresponds to case (ii)
and the top-bottom mirror image of Figure 12 corresponds to case (iii).
(Note that in Figure 12, the two loose ends just to the left of the
generator to be removed are necessarily connected to the south face, by
straight lines.)

The converse statement, that removal of the barred letter in each case
produces an extra loop, is immediate from the diagrams shown.
\qed\enddemo

\remark{Remark 3.1.11}
Consider the special case of Proposition 3.1.9 where the letter $b_p$
to be deleted satisfies $p \in \{1, 2\}$.  Here, there are two
possible cases: (a) that the letter $b_p$ is internal or (b) that we
have deleted an occurrence of $b_2$ which, after commutations, can be
made to appear as the overscored letter in a subsequence of the form 
$b_3 b_1 b_2 b_1 \overline{b_2} b_3$ or $b_3 \overline{b_2} 
b_1 b_2 b_1 b_3$.
In case (b), we say the corresponding occurrence
of $s_2$ is {\it bad}.  Note that all bad occurrences of $s_2$ are
also lateral.
\endremark

\subhead 3.2 The $f$-basis \endsubhead

The aim of \S3.2 is to define a basis of polynomials in the generators
$b_i$ which will eventually be seen to equal the canonical basis in
type $H$.

\proclaim{Lemma 3.2.1}
In a reduced expression $w = s_{i_1} s_{i_2} \cdots
s_{i_l}$ for $w \in W_c$, a bilateral letter in the
expression can only be an occurrence of $s_1$.
\endproclaim

\demo{Proof}
Let $s = s_{i_k}$ be a bilateral letter in the reduced expression.  Since
$s$ is lateral, we can apply commutations so that it lies in a
subsequence $sts$ where $t$ is internal.  Since $w \in W_c$, we must
have $\{s, t\} = \{s_1, s_2\}$, so it remains to show that $s \ne
s_2$.

A bilateral occurrence of $s_2$ would have to occur in a subsequence
$s_1 w' s_2 w'' s_1$, where there are exactly two occurrences of
$s_1$, both internal, and $w'$ and $w''$ are (reduced) words in the
generators that commute with $s_2$.
Since the expression is reduced, $w'$ and $w''$ cannot contain
occurrences of $s_2$.  It follows that $w'$ and $w''$ commute with
$s_1$, which implies that the occurrence of $s_2$ is internal,
contrary to the hypothesis that it is bilateral.
\qed\enddemo

We now state some results on the parsing of reduced expressions in
type $H$ for later usage.

\proclaim{Lemma 3.2.2}
Let $w \in W_c$, and let $s_{i_1} s_{i_2} \cdots s_{i_l}$ be an
arbitrary (but fixed) reduced expression for $w$.

\item{\rm (i)}
{Suppose the expression can be parsed as $w = w_1 s w_2 t w_3 s w_4
t w_5$, where $\{s, t\} = \{s_1, s_2\}$.  Suppose further that $w$ also has
reduced expressions of the forms $$w = w_1 w_2 sts w_3 w_4 t w_5$$ and
$$w = w_1 s w_2 w_3 tst w_4 w_5.$$  Then there is a reduced expression for $w$
for which all four occurrences of $s$ and $t$ occur consecutively.}
\item{\rm (ii)}
{It is impossible to have three consecutive
internal letters in the reduced expression.}
\endproclaim

\demo{Proof}
We first prove (i).  The second and third hypotheses respectively guarantee
that $s$ commutes with
every member of $c(w_2) \cup c(w_3)$ and that $t$ commutes with every
member of $c(w_3) \cup c(w_4)$.  We can therefore transform the
original reduced expression as follows: $$\eqalign{
w &= w_1 s w_2 t w_3 s w_4 t w_5\cr
&= w_1 w_2 st w_3 st w_4 w_5\cr
&= w_1 w_2 w_3 stst w_4 w_5.\cr
}$$  This completes the proof of (i).

To prove (ii), note first that three consecutive internal letters
would have to be of the
form $sts$, where $\{s, t\} = \{s_1, s_2\}$.  Since both occurrences of
$s$ are internal, the $sts$ subsequence must occur in a longer
subsequence of the form $t w' sts w'' t$, where $w'$ and $w''$ are
words in the generators that commute with $t$.  By applying
commutations, we can form a subsequence of the form $tstst$, contrary
to $w \in W_c$.
\qed\enddemo

\definition{Definition 3.2.3}
Let $w \in W_c$.  Let $R$ be the set of internal letters, $s$, of
$w$ that have the property that, after suitable commutations, they 
occur in a subsequence $tst$ where the rightmost occurrence of $t$ is 
bilateral.

Fix a reduced expression for $w \in W_c$.  We say this expression 
is {\it right justified} if (without applying commutations):
\item{(i)}{for all letters $s$ in the set $R$, the letter immediately
to the left is either lateral or internal and}
\item{(ii)}{for all other internal letters $s$, both the letter
immediately to the left and the letter immediately to the right are
either lateral or internal.}
\enddefinition

\example{Example 3.2.4}
Let $w = s_1 s_2 s_3 s_1 s_2 s_1 s_2 \in W(H_3)$ (compare with Example 3.1.6).
The set $R$ has one element, namely the leftmost occurrence of $s_2$.

The reduced expression $s_1 s_2 s_3 s_1 s_2 s_1 s_2$ is right justified, but 
the reduced expression $w = s_1 s_2 s_1 s_3 s_2 s_1 s_2$ is not: the
occurrence of $s_2$ in third from right place fails condition (ii).
\endexample

Right justified reduced expressions for a given $w \in W_c$ need not
be unique, but they always exist, as we see below.

\proclaim{Lemma 3.2.5}
If $w \in W_c$, there exists a right justified reduced expression for
$w$ in which all internal and lateral letters occur in (maximal) contiguous 
subsequences of the form 
(1) $s_1 \overline{s_2}$, (2) $s_1 \overline{s_2} s_1$,
(3) $s_2 \overline{s_1} s_2$,
(4) $s_1 \overline{s_2} \overline{s_1} s_2$,
(5) $s_2 \overline{s_1} \overline{s_2}$,
or (6) $s_2 \overline{s_1} \overline{s_2} s_1$.  (The internal letters
are the barred ones; all other letters are lateral.)
\endproclaim

\demo{Proof}
We describe a simple algorithm to find such a reduced expression.  Start
with any reduced expression and work out which letters are in the set
$R$.  Apply commutations until each internal letter not in the set $R$
occurs as the middle term in a subsequence $sts$.  This produces an
expression of the desired form: note that all the bilateral letters
(which are occurrences of $s_1$ by Lemma 3.2.1) have been commuted to
the right until they are adjacent to internal occurrences of $s_2$.

The assertion about the subsequences comes from the construction and from Lemma
3.2.2.  Part (i) of that lemma explains why it is possible to form
subsequences of types (4), (5) and (6) in the statement, and part (ii)
explains why no longer kinds of subsequence can occur.
\qed\enddemo

We can now define a polynomial over $\zed$ in the generators $b_i$
using the preceding results.

\definition{Definition 3.2.6}
Let $w \in W_c$, and let $s_{i_1} s_{i_2} \cdots s_{i_l}$ be a right
justified reduced expression for $w$.

\vfill\eject

We obtain the polynomial $f_w \in TL(H_{n-1})$ by substituting the
letters in $w$ as follows.  Each of the six subsequences described in Lemma
3.2.5 is replaced as follows: $$\eqalign{
s_1 \overline{s_2} & \ra b_1 b_2 - 1, \cr
s_1 \overline{s_2} s_1 & \ra b_1 b_2 b_1 - b_1, \cr
s_2 \overline{s_1} s_2 & \ra b_2 b_1 b_2 - b_2, \cr
s_1 \overline{s_2} \overline{s_1} s_2 & \ra b_1 b_2 b_1 b_2 - 2b_1 b_2, \cr
s_2 \overline{s_1} \overline{s_2} & \ra b_2 b_1 b_2 - 2b_2, \cr
s_2 \overline{s_1} \overline{s_2} s_1 & \ra b_2 b_1 b_2 b_1 - 2b_2 b_1. \cr
}$$  All other $s_i$ occurring are replaced by the corresponding 
generator $b_i$.  

Note that any bilateral occurrence of $s_1$ must
occur at the beginning of one of the contiguous subsequences listed in
Lemma 3.2.5.  We call such a subsequence {\it distinguished} if it
begins with a bilateral occurrence of $s_1$.  We call a factor of $f_w$
obtained from a distinguished subsequence by the above replacement
procedure a {\it distinguished factor}.
\enddefinition

\example{Example 3.2.7}
Take the right justified reduced expression in Example 3.2.4.  The
corresponding polynomial is $$
f_w = (b_1 b_2 - 1)b_3(b_1b_2b_1b_2 - 2b_1b_2)
.$$  The factor $(b_1b_2b_1b_2 - 2b_1b_2)$ is distinguished because it
arises from a subsequence $s_1 \overline{s_2} \overline{s_1} s_2$
in which the first occurrence of $s_1$ is bilateral (see Example
3.1.6).  The factors
$(b_1 b_2 - 1)$ and $(b_3)$ are not distinguished.
\endexample

\proclaim{Lemma 3.2.8}
\item{\rm (i)}{The elements $\{f_w: w \in W_c\}$ are well-defined and 
form an $\A$-basis of $TL(H_{n-1})$.}
\item{\rm (ii)}{For each $w \in W_c$ we have $f_w \in \L_H$.}
\endproclaim

\demo{Proof}
We tackle (i) first.
The fact that the $f_w$ are well-defined follows from Lemma 3.2.5: all
internal and lateral elements end up in contiguous subsequences
independent of the original reduced expression, and these
are dealt with by Definition 3.2.6.  All other letters go
through a simple replacement procedure.

To see that the set is a basis, expand each $f_w$ in terms of the
monomial basis.  The term $b_w$ occurs with coefficient $1$, and all
other terms are sums of shorter monomials.  Since an arbitrary
monomial is (by using the relations in \S2.1) a linear combination of
basis monomials of equal or shorter length, we see that the set is a basis
as claimed.

We now turn to (ii).  Let $w \in W_c$.  
It is clear from the definition of $f_w$ that it is a $\zed$-linear
combination of monomials in the $b_i$, the longest of which is $b_w$
and the others of which are obtained from $b_w$ by deleting certain
letters.  From the diagram calculus, we see that the operations of
(a) replacing a sequence $b_s b_t b_s$ 
by $b_s$ and (b) replacing a sequence $b_s b_t b_s b_t$ by $b_s b_t$
(where $\{s, t\} = \{s_1, s_2\}$ in both cases) do not alter the number of
loops.  Since $b_w \in \L_H$ by Lemma 3.1.1, it follows that all the
shorter monomials also lie in $\L_H$, and thus that $f_w \in \L_H$, as
required.
\qed\enddemo

\subhead 3.3 Treatment of internal letters \endsubhead

In \S3.3 we show that for each element $f_w$ ($w \in W_c$), there is 
a certain monomial $\wh{f}_w$ in the generators $b_i$ and $\te_i$
which projects via $\pi_H$ to the same element of the lattice $\L_H$.
This is one of three main steps to prove that $\pi_H(f_w) = \pi_H(\te_w)$.

\definition{Definition 3.3.1}
Let $b_w$ ($w \in W_c$) be an element in the monomial basis.  The
{\it expanded form} $b'_w$ of $b_w$ is the monomial 
obtained by replacing each bilateral
occurrence of $b_1$ (see Lemma 3.2.1) by $b_1 b_1$.

We say a monomial $b$ in the generators $b_i$ is $1$-commutative if it
is equal as a monomial, after applying commutations, to $b'_w$ for
some $w \in W_c$.
\enddefinition

\example{Example 3.3.2}
Let $w$ be as in Example 3.2.4.  Then $b = b_1 b_2 b_1 b_3
b_1 b_2 b_1 b_2$ is $1$-commutative, since it is equal to $b_1 b_2 b_3 b_1^2
b_2 b_1 b_2$.
\endexample

The next result
is a version of Proposition 3.1.9 for the expanded forms of monomials.

\proclaim{Lemma 3.3.3}
Let $b = b'_w \in v^k\L_H$ be an expanded monomial.  Let $b_p$ be a
letter in the monomial $b$, and let $\wh{b}$ be the monomial
obtained from $b$ by deleting the letter $b_p$.  Then $\wh{b} \in
v^k \L_H$ unless $b_p$ corresponds to an internal or critical letter of
the basis monomial, $b_w$.
\endproclaim

\demo{Proof}
There are two cases to consider.  In the first case, the deleted letter is
one of the doubled occurrences of $b_1$ arising from a bilateral
occurrence of $b_1$ in $b_w$.  In this case, removing the letter
divides the monomial by $\d$.  Since multiplication by $\d$ increases
the degree of a polynomial in $v$ by 1, we see that deleting this
occurrence of $b_1$ produces a monomial in $v^{k-1} \L_H$, from which
the claim follows.

In the second case, the letter to be deleted is not bilateral.  We
proceed by first applying commutations and the relation $b_1^2 = \d
b_1$, obtaining $b'_w = \d^l b_w$ for some $l$.  The result then
follows from Proposition 3.1.9.
\qed\enddemo

We can modify the $f$-basis in a similar way.

\definition{Definition 3.3.4}
Let $f'_w \in TL(H_{n-1})$ be the element obtained by inserting an
extra occurrence of $b_1$ immediately to the left of each
distinguished factor (see Definition 3.2.6).
\enddefinition

\example{Example 3.3.5}
Let $f_w$ be as in Example 3.2.7.  The underlying word $w$ is given by
$s_1s_2s_3s_1s_2s_1s_2$, and as we can see from Example 3.1.6, the
fourth letter of this word is the only bilateral letter.  Recall from
Example 3.2.7 that the only distinguished factor in $f_w$ is the
rightmost factor in the expression $$
f_w = (b_1 b_2 - 1)(b_3)(b_1 b_2 b_1 b_2 - 2 b_1 b_2)
.$$  Inserting an occurrence of $b_1$ to the left of this factor now
gives $$
f'_w = (b_1b_2 - 1) b_3 b_1 (b_1 b_2 b_1 b_2 - 2 b_1 b_2)
$$ or $$
f'_w = (b_1b_2 - 1)b_1 b_3 (b_1 b_2 b_1 - 2b_1)b_2
.$$
\endexample

\proclaim{Lemma 3.3.6}
We have the following identities:
\item{\rm(i)}{$b_1 b_2 b_1 - b_1 = 
\te_1 \te_2 \te_1
+ v^{-1} \te_2 \te_1
+ v^{-1} v^{-1} \te_1
+ \te_1 \te_2 v^{-1}
+ v^{-1} \te_2 v^{-1}
+ v^{-1} v^{-1} v^{-1}$.}
\item{\rm(ii)}{$b_1 b_2 b_1 b_2 - 2b_1 b_2 = 
(\te_1 \te_2 \te_1 \te_2
+ \te_1 \te_2 \te_1 v^{-1})
+ v^{-1} \te_2 \te_1 \te_2
+ v^{-1} v^{-1} \te_1 \te_2
+ v^{-1} v^{-1} v^{-1} \te_2
+ v^{-1} \te_2 \te_1 v^{-1}
+ v^{-1} v^{-1} \te_1 v^{-1}
+ v^{-1} v^{-1} v^{-1} v^{-1}$.}

There are two other similar identities obtained
by exchanging the roles of $1$ and $2$ above.
\endproclaim

\demo{Note}
{The reason the result has been expressed in such an inconcise
way will become clear in the proof of Lemma 3.3.9.}
\enddemo

\demo{Proof}
This is a routine calculation.
\qed\enddemo

\definition{Definition 3.3.7}
Let $w \in W_c$.  The element $\wh{f}_w \in TL(H_{n-1})$ is that
obtained by taking the corresponding monomial basis element $$
b_w = b_{i_1} b_{i_2} \cdots b_{i_r}
$$ and substituting $\te_i$ for $b_i$ whenever $b_i$ is an internal or
lateral letter that is not a bad occurrence of $b_2$ (see Remark 3.1.11).

We define an expanded form $\wh{f}'_w$ by
replacing each occurrence of $\te_1$ in $\wh{f}_w$ that corresponds to
a bilateral occurrence of $s_1$ in $w$ by $\te_1 \te_1$.
\enddefinition

\example{Example 3.3.8}
Let $w = s_1 s_2 s_3 s_1 s_2 s_1 s_2 s_3\in W(H_3)$.  In this case, the
rightmost occurrence of $s_2$ is bad, and we have $$
\wh{f}_w = \te_1 \te_2 b_3 \te_1 \te_2 \te_1 b_2 b_3
$$ and $$\wh{f}'_w
= \te_1 \te_2 b_3 \te_1 \te_1 \te_2 \te_1 b_2 b_3
.$$  It is more convenient for later purposes to move the new
occurrence of $\te_1$ as far to the left as possible; this gives $$
\wh{f}'_w 
= \te_1 \te_2 \te_1 b_3 \te_1 \te_2 \te_1 b_2 b_3
.$$
\endexample

\proclaim{Lemma 3.3.9}
For each $w \in W_c$, let $\k = \k(w)$ be the number of bilateral
occurrences of $s_1$ in $w$.  Then 
$\pi_H(v^{-\k}f'_w) = \pi_H(v^{-\k}\wh{f}'_w)$, where $\pi_H$ is as in
Definition 3.1.2.
\endproclaim

\demo{Proof}
Consider $f'_w$, and replace all the sections corresponding (via
Definition 3.2.6) to internal and lateral letters using the identities
in Lemma 3.3.6.  This shows that $f'_w$ is equal to $\wh{f}'_w$ plus
lower monomial terms (in $\te_i$, $b_i$ and $v^{-1}$)
which differ from $f'_w$ in that certain of the $\te_i$
have been replaced by $v^{-1}$.  Note that each term containing a $v^{-1}$ is
missing at least one lateral $\te_i$.

We now expand each of the lower monomial terms above to a polynomial in
the $b_i$ by using the identity $\te_i = b_i - v^{-1}$.  This
expression expands to a sum of monomials in the $b_i$, each of which
can be obtained from the expanded form $b'_w$ by (possibly) applying a
sign change and then deleting some of the letters and replacing them
with instances of $v^{-1}$.  As in the previous paragraph,
one of the missing letters is guaranteed to be lateral.  Note that
$b'_w = \d^\k b_w \in v^\k \L_H$.  

If the missing lateral letter
is not a bad occurrence of $b_2$ then it is not internal or critical.
In this case, we use Lemma 3.3.3 (applied to the missing lateral letter) 
and apply Lemma 3.1.3 repeatedly (if necessary) to show 
that all the lower monomial terms lie in $v^{\k - 1}\L_H$.
It follows that after multiplication by $v^{-\k}$ and applying $\pi_H$, the
lower terms go to zero, which completes the argument.

If the missing lateral letter is a bad occurrence of $b_2$ then we
proceed in the same way, but we are left with an extra term which will
be $\te_1 \te_2 \te_1 v^{-1}$ or its left-right mirror image.  The
terms in Lemma 3.3.6 (ii) will all go to zero as before, except the
two in parentheses; these can be combined to form $\te_1 \te_2 \te_1
b_2$ or its left-right mirror image; this is the term we require.
\qed\enddemo

\proclaim{Lemma 3.3.10}
For each $w \in W_c$, $\pi_H(f_w) = \pi_H(\wh{f}_w)$.
\endproclaim

\demo{Proof}
Let $\k = \k(w)$ as in Lemma 3.3.9.  It is clear that
$\pi_H(v^{-\k}f'_w) = \pi_H(f_w)$, so by Lemma 3.3.9, it is enough to
prove that $\pi_H(v^{-\k}\wh{f}'_w) = \pi_H(\wh{f}_w)$.

Observe that $v^{-1} (\te_1)^2 = v^{-1} + (1 - v^{-2}) \te_1$.  We
need to apply
this relation $\k$ times to $\wh{f}'_w$.  Arguing as in the proof of
Lemma 3.3.9, we find that the $v^{-1}$ which appears
on the right hand side can be ignored (it corresponds to the deletion
of a non-critical lateral letter and replacement by $v^{-1}$), as can 
the term $v^{-2} \te_1$.  The claim follows.
\qed\enddemo

\subhead 3.4 Agreement of the $f$-basis and the canonical basis \endsubhead

%In \S3.4 we show how more of the $b_i$ in the monomial $\wh{f}_w$ may
%be converted into $\te_i$ without changing the image of the monomial
%under $\pi_H$.  This is the second of the three main steps to prove that
%$\pi_H(f_w) = \pi_H(\te_w)$. 

In \S3.4, we will often implicitly use the fact that $b_i$ and $\te_j$
commute if and only if $b_i$ and $b_j$ commute if and only if $\te_i$
and $\te_j$ commute.

\proclaim{Lemma 3.4.1}
Let $w \in W_c$, and let $\wh{f}_w$ be as in Definition 3.3.7.
Let $f$ be a monomial in the $b_i$ and $\te_j$ obtained from
$\wh{f}_w$ by changing one critical occurrence of $b_p$ to $\te_p$.
Then $\pi_H(f) = \pi_H(\wh{f}_w)$.
\endproclaim

\demo{Proof}
It is convenient to split the proof into two cases: $p = 2$ (the case of the
bad occurrences of $2$) and $p =
2k > 2$.  There are no other possibilities (see Proposition 3.1.9).

Suppose $p = 2$ and consider the difference $d = \wh{f}_w - f$.  This is
obtained by replacing the deleted occurrence of $b_p$ in $\wh{f}_w$ by 
$v^{-1}$.  By Remark 3.1.11, we may apply commutations to 
$d$ so that it is equal to $x b_3 \te_1 \te_2 \te_1
v^{-1} b_3 x'$ or its mirror image, where $x$ and $x'$ are monomials
in the $b_i$ and $\te_i$.  We
treat only the first case; the other is the same.  Now $$
\te_1 \te_2 \te_1 = (b_{121} - b_1) - v^{-1} b_2 b_1 - v^{-1} b_1 b_2
+ v^{-2} b_1 + v^{-2} b_2 - v^{-3}
$$ which gives $$\eqalign{
x b_3 \te_1 \te_2 \te_1 v^{-1} b_3 x' =& 
x v^{-1} b_3(b_{121} - b_1)b_3 x' \cr
&- x b_3 v^{-1} b_2 b_1 v^{-1} b_3 x'
- x b_3 v^{-2} b_1 b_2 b_3 x' \cr 
&+ x b_3 v^{-2} b_1 v^{-1}b_3 x' + 
x b_3 v^{-3} b_2 b_3 x' 
- x b_3 v^{-4} b_3 x'.\cr
}$$  The first term in this sum is identically zero, and the others
project under $\pi_H$ to zero by the argument of Lemma 3.3.9.  This
completes the proof when $p = 2$.

Suppose now that $p = 2k > 2$.  In this case, we may apply
commutations to $\wh{f}_w$ and Proposition 3.1.9 (part (i), (ii) or (iii))
to obtain an expression $$
x b_3 b_5 \cdots b_{2k-1} \te_1 \te_2 \te_1 b_4 b_6 \cdots
b_{2k} b_3 b_5 \cdots b_{2k-1} x'
$$ if we are in case (i) of the proposition, $$
x b_3 b_5 \cdots b_{2k+1} \te_1 \te_2 \te_1 b_4 b_6 \cdots
b_{2k} b_3 b_5 \cdots b_{2k-1} x'
$$ in case (ii) and $$
x b_3 b_5 \cdots b_{2k-1} \te_1 \te_2 \te_1 b_4 b_6 \cdots
b_{2k} b_3 b_5 \cdots b_{2k+1} x'
$$ in case (iii).  In all cases, the positions of $b_3$ ensure that
the $\te_1$ occurring are lateral occurrences and the $\te_2$ is internal.
We may now invoke the argument used to deal with the $p = 2$ case to prove
that we may replace the string $\te_1 \te_2 \te_1$ by $(b_{121} -
b_1)$ without changing the result after projection by $\pi_H$.  Removing
the generator $b_p$ from the expression and replacing it by $v^{-1}$
to produce $d$ now yields, in each case, an expression with $$
b_3 b_5 \cdots b_{2l + 1} (b_{121} - b_1) b_4 b_6 \cdots b_{2l} b_3
b_5 \cdots b_{2l + 1}
$$ somewhere in the middle, where we have $l = k-1$ in case (i) and $l = k$ in
cases (ii) and (iii).  This expression is identically zero; this can
be seen from the diagram calculus.  This proves that $\pi_H(d) = 0$ and 
completes the proof in the case $p = 2k > 2$.
\qed\enddemo

\proclaim{Corollary 3.4.2}
Let $w \in W_c$, and let $\wh{f}_w$ be as in Definition 3.3.7.
Let $\ti{f}_w$ be the monomial in the $b_i$ and $\te_j$ obtained from
$\wh{f}_w$ by changing all the critical occurrences of $b_p$ to $\te_p$.
Then $\pi_H(\ti{f}_w) = \pi_H(\wh{f}_w)$.
\endproclaim

\demo{Proof}
The proof is by repeated applications of Lemma 3.4.1, starting by
replacing the critical occurrences of $b_p$ for the largest possible
$p$.  This works because $p$ must be even and, in the proof of Lemma 3.4.1,
removal of $b_p$ only involves generators $b_q$ and $\te_q$ with 
$q < p+2$.
\qed\enddemo

We can now prove the main result of \S3.

\proclaim{Theorem 3.4.3}
The basis $\{f_w : w \in W_c\}$ is the canonical basis of
$TL(H_{n-1})$ in the sense of Theorem 1.2.5.
\endproclaim

\demo{Proof}
Since $f_w$ is a $\zed$-linear combination of monomials in the $b_i$,
it follows that they are fixed by $\bar{\ }$.  It remains to show that
$\pi(f_w) = \pi(\te_w)$ for all $w \in W_c$, where $\pi$ is the
projection defined in \S1.2.

By Lemma 3.3.10, $\pi_H(f_w) = \pi_H(\wh{f}_w)$, and by Corollary 3.4.2,
$\pi_H(\wh{f}_w) = \pi_H(\ti{f}_w)$.  It remains to show that
$\pi_H(\ti{f}_w) = \pi_H(\te_w)$.  This is achieved by the argument in
Lemma 3.3.9: we apply Lemma 3.3.3 and then apply Lemma 3.1.3
repeatedly.  This is valid since all the problem cases---the internal and
critical letters described in Proposition 3.1.9---have been dealt with.

It follows that $\te_w = \sum_{x \in W_c} c_x f_x$, where $c_x \in
v^{-1} \A^-$ for $x \ne w$ and $c_w - 1 \in v^{-1} \A^-$.  A
standard argument (as in the proof of \cite{{\bf 6}, Theorem 3.6}) shows
that we also have $f_w = \sum_{x \in W_c} c'_x \te_x$, where $c'_x$
satisfies the same properties as $c_x$ above.  This shows that
$\pi(f_w) = \pi(\te_w)$, which completes the proof.
\qed\enddemo

\subhead 4. Properties of the basis in type $H$ \endsubhead

The aim of \S4 is to show that the $f$-basis constructed in \S3 is
equal to the diagram basis described in \S2.  This will prove Theorem 2.1.3.

\subhead 4.1 Positivity properties for the $f$-basis \endsubhead

The aim of \S4.1 is to establish a positivity property for the structure
constants of $TL(H_{n-1})$ with respect to the $f$-basis.
As in \S3, all computations in \S4 take place in $TL(H_{n-1})$ over the ring
$\A$ unless otherwise stated.

\proclaim{Lemma 4.1.1} Let $w \in W_c$.
\item{\rm (i)}{Suppose $s \in S$ is such that $ws \notin W_c$.  Then
there exists a unique $s' \in S$ such that any reduced expression can
be parsed in one of the following two ways:
\item{\rm (a)}{$w=w_1sw_2s'w_3$, where $ss'$ has order 
$3$, and $s$ commutes with every member of 
$c(w_2)\cup c(w_3)$;}
\item{\rm (b)}{$w=w_1sw_2s'w_3sw_4s'w_5$, where $\{s, s'\} = \{s_1,
s_2\}$, $s$ commutes with every member of 
$c(w_2)\cup c(w_3) \cup c(w_4) \cup c(w_5)$, and $s'$ commutes 
with every member of $c(w_3)\cup c(w_4)$.}
}
\item{\rm (ii)}{Suppose $s \in S$ is such that $ws \in W_c$ and the
occurrence of $s$ shown is lateral.  Let $s'$ be such that $\{s, s'\}
= \{s_1, s_2\}$.  Then any reduced expression can be parsed in one of the
following two ways:
\item{\rm (a)}{$w=w_1\overline{s}w_2s'w_3$, where $s$ commutes with 
every member of $c(w_2)\cup c(w_3)$ and the overscored occurrence of
$s$ is not internal;}
\item{\rm (b)}{$w=w_1 \overline{s'} w_2 s w_3s'w_4$, where $s$ commutes with 
every member of $c(w_3)\cup c(w_4)$, $s'$ commutes with every member
of $c(w_2)\cup c(w_3)$ and the overscored occurrence of
$s'$ is not internal.}
}
\endproclaim

\demo{Proof}
The proof of (i) follows exactly the same principles as the
corresponding argument in type $B$.  Since the latter argument is presented in
full detail in \cite{{\bf 7}, Lemma 2.1.2}, we omit the proof.

The proof of (ii) follows the same principles.
We may apply commutations to the element $ws$ to form a maximal
contiguous sequence $ss's$ or $s'ss's$ in which the original
occurrence of $s$ appears on the right; no other cases are possible by 
Lemma 3.2.2 (ii).  The first case is described by (a) and the second
by (b).  (Compare also with Lemma 3.2.5.)
\qed\enddemo

Lemma 4.1.1 can be refined if we assume the reduced expressions
occurring are right justified.

\proclaim{Corollary 4.1.2}
If all reduced expressions occurring in Lemma 4.1.1 are right
justified, the situations in parts (i)(b), (ii)(a) and (ii)(b) of that
result may be respectively simplified to
\item{\rm (i)}{$w=w_1ss'ss'w_2$, where $\{s, s'\} = \{1,
2\}$ and $s$ commutes with every member of $c(w_2)$;
}
\item{\rm (ii)}{$w=w_1\overline{s}s'w_2$, where $s$ commutes with 
every member of $c(w_2)$ and the overscored occurrence of
$s$ is not internal;}
\item{\rm (iii)}{$w=w_1 \overline{s'} s s'w_2$, where $s$ commutes with 
every member of $c(w_2)$  and the overscored occurrence of
$s'$ is not internal.}

In all three cases, $s$ and $s'$ do not lie in $c(w_2)$.
\endproclaim

\demo{Proof}
This follows quickly from the definition of ``right justified''
and Lemma 4.1.1.
\qed\enddemo

\definition{Definition 4.1.3}
We set $\A^{\geq 0} := \enn[v, v^{-1}]$, where the natural numbers
$\enn$ are taken to include zero.

If $a, b \in W(X)$ are comparable in the Bruhat--Chevalley order, we define
\newline $\max(a, b)$ to be $a$ if $a \geq b$ and $b$ if $b \geq a$.
\enddefinition

\proclaim{Lemma 4.1.4}
Let $w \in W_c$.  Then the structure constants $a_x$ in the 
expression $$
f_w b_i = \sum_{x \in W_c} a_x f_x
$$ satisfy $a_x \in \A^{\geq 0}$.  Furthermore, for all $a_x \ne 0$, we
have $x \leq \max(w, ws_i)$ and 
$x s_i < x$, where $<$ is the Bruhat--Chevalley order.
\endproclaim

\demo{Proof}
We claim that the statement is true for $n = 3$, \idest $TL(H_2)$.
In this case it may be verified that the $f$-basis is equal to the
diagram basis and the other assertions follow from a case-by-case check.

For the case of general $n$, we proceed by induction on $\ell(w)$.  
If $\ell(w) < 2$, the
statement is obvious.  Now consider $f_w b_i$ where $\ell(w) = k > 1$ and
the statement is known to be true for $\ell(w) < k$.

There are three cases: the first is $w s_i < w$;
the second is $w s_i > w$ and $w s_i \in W_c$;
the third is $w s_i > w$ but $w s_i \not\in W_c$.

In the first case, $w$ has a reduced expression ending with $s_i$.
It follows from the construction of the $f$-basis in Definition 3.2.6
that $f_w$ is equal to an element of the form $f' b_i$ for some $f'$.
This means $f_w b_i = (v + v^{-1}) f_w$ and the hypotheses are satisfied.

Next, we tackle the second case.  There are two subcases
according as the rightmost occurrence of $s_i$ is lateral or not.  If
the occurrence is not lateral, we find that $f_{ws_i} = f_w b_i$, and
we are done.  If the occurrence is lateral, we are in the situation
described in Lemma 4.1.1 (ii).  We use the simplification provided by
Corollary 4.1.2 and the $n = 3$ case mentioned at the beginning of the
proof and set $s := s_i$.  In the case of Lemma 4.1.1 (ii) (a), $w =
w_1ss'w_2$  and  $$
f_w b_i = f_{w_1} (f_{ss'} b_i) f_{w_2} = f_{w_1} (f_{ss's} + f_s) f_{w_2}
= f_{ws} + (f_{w_1s} f_{w_2})
.$$  The term $f_{w_1s}f_{w_2} = f_{w_1s}b_{w_2}$ is equal by
induction to an $\A^{\geq 0}$-linear combination of basis elements $f_x$
for which $x$ is an ordered subexpression of $w_1sw_2$, and therefore
$x \leq w$.
Since $s$ commutes with all elements in $c(w_2)$, the basis elements
$f_x$ occurring also satisfy $xs < x$.  This completes the first subcase.  The
second subcase, corresponding to Lemma 4.1.1 (ii) (b), is similar but
uses the identity $w = w_1 s'ss' w_2$ and $$
f_w b_i = f_{w_1} (f_{s'ss'} b_i) f_{w_2} = f_{w_1} (f_{s'ss's} +
f_{s's}) f_{w_2} = f_{ws} + (f_{w_1ss'} f_{w_2})
.$$

The third case is rather similar to the second case, except that there
is no term $f_{ws}$.  The identities
to use here are $w = w_1ss'ss'w_2$ and $$
f_w b_i = f_{w_1} (f_{ss'ss'} b_i) f_{w_2} = f_{w_1} f_{ss's} f_{w_2}
= f_{w_1ss's} f_{w_2}
.$$  \qed\enddemo

\proclaim{Corollary 4.1.5}
Let $w \in W_c$.  Then $f_w b_i = (v + v^{-1}) f_w$ if and only if
$ws_i < w$.
\endproclaim

\demo{Proof}
If $ws_i < w$ then the proof of Lemma 4.1.4 shows that $f_w b_i = (v +
v^{-1}) f_w$. For the converse, consider the product $f_w b_i$.  In
the second case of the proof of Lemma 4.1.4, this product contains
$f_{ws}$ with coefficient 1.  In the third case of the proof, the
product is a linear combination of $f_x$ where $x$ is shorter than
$w$.  By the positivity property in Lemma 4.1.4, no cancellation can occur
and it is impossible in these cases for $f_w b_i = (v + v^{-1}) f_w$.  We must
therefore be in case (i) of the proof, which gives $ws_i < w$.
\qed\enddemo

\proclaim{Lemma 4.1.6}
Let $w \in W_c$.  Suppose that $ws' < w$ for $s' \in \{s_1, s_2\} =
\{s, s'\}$. 
\item{\rm (i)}{The structure constants $a_x$ in the 
expression $$
f_w (b_s b_{s'} - 1) = \sum_{x \in W_c} a_x f_x
$$ satisfy $a_x \in \A^{\geq 0}$.}
\item{\rm (ii)}{The structure constants $a'_x$ in the 
expression $$
f_w (b_s b_{s'} b_s - 2b_s) = \sum_{x \in W_c} a'_x f_x
$$ satisfy $a'_x \in \A^{\geq 0}$.}
\endproclaim

\demo{Proof}
We proceed by induction on $\ell(w)$, the cases
of length $0$ and $1$ being easy.

Since $w \in W_c$ has a reduced expression ending in $s'$, it has none ending
in $s$, so $ws > w$.  There are three cases to consider: first if
$ws \in W_c$ and the occurrence of $s$ shown is not lateral; second
if $ws \in W_c$ and the occurrence of $s$ shown is lateral and third
if $ws \not\in W_c$.

For the first case, we first observe that $wss'$ and $wss's$
are also in $W_c$.  This is a consequence of Lemma 4.1.1 (ii).
The definition of the $f$-basis shows that $f_{wss'} = f_w (b_s b_{s'}
- 1)$, which proves (i), and that $f_{wss's} = f_w (b_{ss's} -
2b_s)$, which proves (ii).

For the second case, we proceed as in the second case in the proof of
Lemma 4.1.4.  There are four subcases, corresponding to combinations
of either parts (ii) or (iii) of Corollary 4.1.2 and either part (i)
or part (ii) of the statement.  Let us deal with part (i) of the
statement.  Here, we have $w =
w_1\overline{s}s'w_2$ or $w = w_1\overline{s'}ss'w_2$ where 
$s \not\in c(w_2)$ and the overscored letter is not internal.
Since we also assume $ws' < w$,
we must have $s' \not\in c(w_2)$ as well.  We now have $$\eqalign{
f_w (b_s b_{s'} - 1) &= f_{w_1} f_{ss'} f_{w_2} (b_s b_{s'} - 1) =
f_{w_1} f_{ss'} (b_s b_{s'} - 1) f_{w_2}\cr
&= f_{w_1} (f_{ss'ss'} + f_{ss'}) f_{w_2} = f_{wss'} + f_{w_1ss'} f_{w_2}.
}$$  Since $f_{w_2} = b_{w_2}$, we are done by Lemma 4.1.4.  The other
subcase is similar but uses the fact that $f_{s'ss'}(b_sb_{s'} - 1) =
f_{s'ss'} + f_{s'}$.
The proof for part (ii) follows a similar pattern, again relying on
positivity properties in the case $n = 3$: $f_{ss'}(b_s b_{s'} b_s -
2b_s) = f_{ss's}$ and $f_{s'ss'}(b_s b_{s'} b_s - 2b_s) = f_{s's}$.

The third case is reminiscent of the second case.  Here, $w = w_1
\overline{s}s'ss'w_2$; we need the facts that
$f_{ss'ss'}(b_sb_{s'} - 1) = f_{ss'}$ and $f_{ss'ss'}(b_s b_{s'} b_s -
2b_s) = f_s$.
\qed\enddemo

There are also left-handed versions of the results in lemmas 4.1.4 and
4.1.6.  The proofs of these follow by making trivial changes to the
arguments, but we present the statements for completeness in the
following two lemmas.

\proclaim{Lemma 4.1.7}
Let $w \in W_c$.  Then the structure constants $a_x$ in the 
expression $$
b_i f_w = \sum_{x \in W_c} a_x f_x
$$ satisfy $a_x \in \A^{\geq 0}$.  Furthermore, for all $a_x \ne 0$, we
have $x \leq \max(w, s_iw)$ and 
$s_i x < x$, where $<$ is the Bruhat--Chevalley order.
\endproclaim

\proclaim{Lemma 4.1.8}
Let $w \in W_c$.  Suppose that $s'w < w$ for $s' \in \{s_1, s_2\} =
\{s, s'\}$. 
\item{\rm (i)}{The structure constants $a_x$ in the 
expression $$
(b_{s'} b_s - 1) f_w = \sum_{x \in W_c} a_x f_x
$$ satisfy $a_x \in \A^{\geq 0}$.}
\item{\rm (ii)}{The structure constants $a'_x$ in the 
expression $$
(b_s b_{s'} b_s - 2b_s)f_w = \sum_{x \in W_c} a'_x f_x
$$ satisfy $a'_x \in \A^{\geq 0}$.}
\endproclaim

These results have an important consequence which verifies \cite{{\bf 7},
Conjecture 1.2.4} for Coxeter systems of type $H$:

\proclaim{Proposition 4.1.9}
The structure constants for the canonical basis for $TL(H_{n-1})$ 
lie in $\A^{\geq 0}$.
\endproclaim

\demo{Proof}
By the construction of the $f$-basis, we see that $f_w$ is either:
\item{(a)}{of
the form $b_i f_{w'}$ or $f_{w'}b_i$ for some $w' \in W_c$ such that
$\ell(w') = \ell(w) - 1$, or}
\item{(b)}{of the form $(b_{s'} b_s - 1)f_{w'}$ or
$f_{w'} (b_s b_{s'} - 1)$ for
some $w' \in W_c$ such that $s'w' < w'$ (respectively, $w's' < w'$)
and $\ell(w') = \ell(w) - 2$ or}
\item{(c)}{of the form $(b_s b_{s'} b_s - 2b_s)f_{w'}$ or
$f_{w'} (b_s b_{s'} b_s - 2b_s)$ for
some $w' \in W_c$ such that $s'w' < w'$ (respectively, $w's' < w'$) 
and $\ell(w') = \ell(w) - 3$.}

It follows that the multiplication of two basis elements $f_x f_y$ for
$x, y \in W_c$ can be broken down into repeated applications of lemmas
4.1.4, 4.1.6, 4.1.7 and 4.1.8, from which the result follows.
\qed\enddemo

Note that Proposition 4.1.9 is also a consequence of the as yet
unproved Theorem 2.1.3 together with \cite{{\bf 5}, Proposition 4.1.1}.

\subhead 4.2 Comparison with the diagram basis \endsubhead

The diagram basis for $TL(H_{n-1})$ can easily be shown to satisfy
properties analogous to those described in \S4.1 by using results from
\cite{{\bf 5}}.

We denote diagram basis elements by $D_w$ (although we do not specify
a bijection between $W_c$ and the set of diagrams).

The following result, which is analogous to lemmas 4.1.4 and 4.1.7, is
immediate from the diagram calculus.

\proclaim{Lemma 4.2.1}
The structure constants $a_x$ in the expression $$
D_w b_i = \sum_{x \in W_c} a_x D_x
$$ satisfy $a_x \in \A^{\geq 0}$.  Furthermore, for all $a_x \ne 0$, we
have $D_x b_i = (v + v^{-1}) D_x$.

An analogous property holds for $b_i D_w$.
\endproclaim

The next result is the analogue of lemmas 4.1.6 and 4.1.8.  The proof
follows easily from the diagram calculus.

\proclaim{Lemma 4.2.2}
Suppose that $D_w b_{s'} = (v + v^{-1}) D_w$ for $s' \in \{s_1, s_2\} =
\{s, s'\}$. 
\item{\rm (i)}{The structure constants $a_x$ in the 
expression $$
D_w (b_s b_{s'} - 1) = \sum_{x \in W_c} a_x D_x
$$ satisfy $a_x \in \A^{\geq 0}$.}
\item{\rm (ii)}{The structure constants $a'_x$ in the 
expression $$
D_w (b_s b_{s'} b_s - 2b_s) = \sum_{x \in W_c} a'_x D_x
$$ satisfy $a'_x \in \A^{\geq 0}$.}

There are also analogous results for $(b_{s'}b_s - 1)D_w$ and
$(b_sb_{s'}b_s - 2b_s)D_w$.
\endproclaim

The diagram basis has the following key property.

\proclaim{Lemma 4.2.3}
Each diagram basis element $D$ may be obtained from the identity
element $1$ by repeated application of the following six procedures:
\item{\rm(i)}{left multiplication by $b_i$;}
\item{\rm(ii)}{right multiplication by $b_i$;}
\item{\rm(iii)}{left multiplication of an element $D'$ by $b_{s'}b_s -
1$ where $b_{s'}D' = (v+v^{-1})D'$ and $\{s, s'\} = \{s_1, s_2\}$;}
\item{\rm(iv)}{right multiplication of an element $D'$ by $b_sb_{s'} -
1$ where $D'b_{s'} = (v+v^{-1})D'$ and $\{s, s'\} = \{s_1, s_2\}$;}
\item{\rm(v)}{left multiplication of an element $D'$ by $b_s b_{s'}b_s -
2b_s$ where $b_{s'}D' = (v+v^{-1})D'$ and $\{s, s'\} = \{s_1, s_2\}$;}
\item{\rm(vi)}{right multiplication of an element $D'$ by $b_s b_{s'}b_s -
2b_s$ where $D'b_{s'} = (v+v^{-1})D'$ and $\{s, s'\} = \{s_1, s_2\}$;}
\endproclaim

\demo{Proof}
This is implicit in the proof of \cite{{\bf 5}, Proposition 3.2.8}; see
also the other results of \cite{{\bf 5}, \S3.2}.
\qed\enddemo

We are now ready to prove Theorem 2.1.3.

\demo{Proof of Theorem 2.1.3}
Combining Lemma 4.2.3 with 
lemmas 4.1.4, 4.1.6, 4.1.7, 4.1.8 and Corollary 4.1.5 shows that each
diagram basis element $D$ is an $\A^{\geq 0}$-linear combination of
$f$-basis elements.  

Conversely, combining the construction of the
$f$-basis as in the proof of Proposition 4.1.9 with Corollary 4.1.5
and lemmas 4.2.1 and 4.2.2 shows that each $f$-basis element is an
$\A^{\geq 0}$-linear combination of diagram basis elements.

We observe that a square matrix with entries in $\A^{\geq 0}$ whose
inverse has entries in $\A^{\geq 0}$ must be a monomial matrix whose
entries are of the form $v^k$ for some $k \in \zed$.  Lemma 4.2.3 and the
definition of the $f$-basis show that both the $f$-basis and the
diagram basis are invariant under the map $\bar{ }$ of Definition
1.2.4.  It follows that all entries in the transition matrices between
the $f$-basis and the diagram basis must be invariant under $\bar{ }$, and
thus that all entries of the monomial matrix are equal to $1$ (\idest
$k = 0$ above).  The
transition matrices are therefore permutation matrices and the
$f$-basis equals the diagram basis.  The proof is completed by Theorem 3.4.3.
\qed\enddemo

\subhead 5. Combinatorics in type $B$  \endsubhead

In this final section, we show how the arguments of \S3 and \S4 can be 
adapted to work for type $B$, and thus to prove Theorem 2.2.5.  The 
argument is essentially a subset of
the argument for type $H$, but we have chosen to present the case of
type $B$ separately to avoid making the earlier arguments overly complicated.

In \S5, all computations take place in $TL(B_{n-1})$ over the ring
$\A$ unless otherwise stated.

\subhead 5.1 Basic properties of $B$-canonical diagrams \endsubhead

\proclaim{Lemma 5.1.1}
The set $C_n$ of Theorem 2.2.5 
is linearly independent over $\A$, and its structure constants lie
in $\A^{\geq 0}$.
\endproclaim

\demo{Proof}
The fact that the set is linearly independent follows from comparing
definitions 2.2.2 and 2.2.4.

It is obvious from the relations in figures 4 and 5 that the 
structure constants for this set lie in $\A[{1\over2}]$.  

The proof is a case by case check, multiplying elements of $C_n$
of various types (C$1$, C$1'$, C$2$ as in Definition 2.2.4) together.  The two
difficult cases are (a) removal of a (single) loop decorated by a circle and
(b) what to do when a diagram emerges with a (single instance of a) 
circular decoration where
a square is required.  In each of these cases, there is a spare factor
of $2$ available; this deals with case (a) immediately.  For case (b),
we apply the first relation in Figure 5 to reexpress the diagram as a
linear combination of two elements of $C_n$, each with coefficient $1$.
\qed\enddemo

The next problem to resolve is the issue of which $\A$-form of the
algebra $TL(B_{n-1})$ to use.  (The one given in \cite{{\bf 4}} is the
wrong one for our purposes.)

We start by defining an injective map from the set of $B$-canonical
diagrams (see Definition 2.2.4) to the set of $H$-admissible diagrams
(see Definition 2.1.1).

\definition{Definition 5.1.2}
Let $D$ be a $B$-canonical diagram, and let $\l D \in C_n$ (so that
$\l = 1$ or $2$).  Then the linear map $\iota : TL(B_{n-1}) \ra
TL(H_{n-1})$ is defined to take $\l D$ to the $H$-admissible diagram obtained 
by replacing all the decorations in $D$ (whether
circular or square) by circular decorations.
\enddefinition

\proclaim{Lemma 5.1.3}
The map $\iota$ is well-defined and injective.
\endproclaim

\demo{Proof}
This follows from examination of definitions 2.1.1 and 2.2.4.
\qed\enddemo

The next part of the argument is an adaptation of the argument in
\cite{{\bf 5}, \S3.2} to type $B$.
The proofs of the following four lemmas, in which $D \in C_n$,
are immediate.

\proclaim{Lemma 5.1.4}
Assume $D$ has a propagating edge, $E$, connecting node $p_1$ in the north
face to node $p_2$ in the south face.

If nodes $p_1 + 1$ and $p_1 + 2$ in the north face
are connected by a (necessarily
undecorated) edge $E'$, then $b_{p_1} D$ is the element of $C_n$
obtained by removing $E'$, 
disconnecting $E$ from the north face and reconnecting it to
node $p_1 + 2$ in the north face,
and installing a new undecorated edge between
points $p_1$ and $p_1 + 1$ in the north face.  The edge corresponding
to $E$ retains
its original decoration status.

Furthermore, $\i(b_{p_1}D) = b_{p_1}\i(D)$ in $TL(H_{n-1})$.
\endproclaim

\proclaim{Lemma 5.1.5}
Assume that $i > 1$, and that in the north face of $D$, nodes
$i$ and $i+1$ are connected by a (square-)decorated edge $e_1$,
and nodes $i+2$ and $i+3$ are connected
by an undecorated edge, $e_2$.  Then $b_i b_{i+1}
D$ is the element of $C_n$ obtained from $D$ by exchanging $e_1$
and $e_2$.  This procedure has an inverse, since $D = b_{i+2} b_{i+1}
b_i b_{i+1} D$.

Furthermore, $\i(b_i b_{i+1}D) = b_i b_{i+1}\i(D)$ in $TL(H_{n-1})$.
\endproclaim

\proclaim{Lemma 5.1.6}
Assume that in the north face of $D$, nodes 1 and 2 are connected by
a decorated edge (necessarily decorated by a circle), and nodes 3 and
4 are connected by an undecorated
edge.  Then $(b_1b_2 - 1) D \in C_n$ is the element obtained
from $D$ by decorating the edge connecting nodes 3 and 4 with a square.

Furthermore, $\i((b_1b_2 - 1)D) = (b_1b_2 - 1)\i(D)$ in $TL(H_{n-1})$.
\endproclaim

\proclaim{Lemma 5.1.7}
Assume that in the north face of $D$, nodes $i$ and $i+1$ are
connected by an undecorated edge, $e_1$, and nodes $j < i$ and $k > i+1$ are
connected by an edge, $e_2$.  Assume also that $j$ and $k$ are chosen
such that $|k - j|$ is minimal.  Then $D$ is of the form $b_i D'$, where
$D'$ is an element of $C_n$ which 
is the same as $D$ except as regards the edges connected to nodes $j,
i, i+1, k$ in the north face.  Nodes $j$ and $i$ in $D'$ are connected
to each other by an edge with the same decoration as $e_2$, and nodes
$i+1$ and $k$ are connected to each other by an undecorated edge.

Furthermore, $\i(b_i D') = b_i\i(D')$ in $TL(H_{n-1})$.
\endproclaim

These results may be visualized as in Figure 13, in which an unshaded
circular decoration denotes an optional decoration of indeterminate
type, thus depicting the concept of ``original decoration status''
in Lemma 5.1.4.  There are also
right-handed versions of the results, which correspond to top-bottom
reflected versions of Figure 13.

\vfill\eject

\topcaption{Figure 13} Lemmas 5.1.4--5.1.7 \endcaption
\centerline{
\hbox to 4.888in{
\vbox to 2.708in{\vfill
        \includegraphics{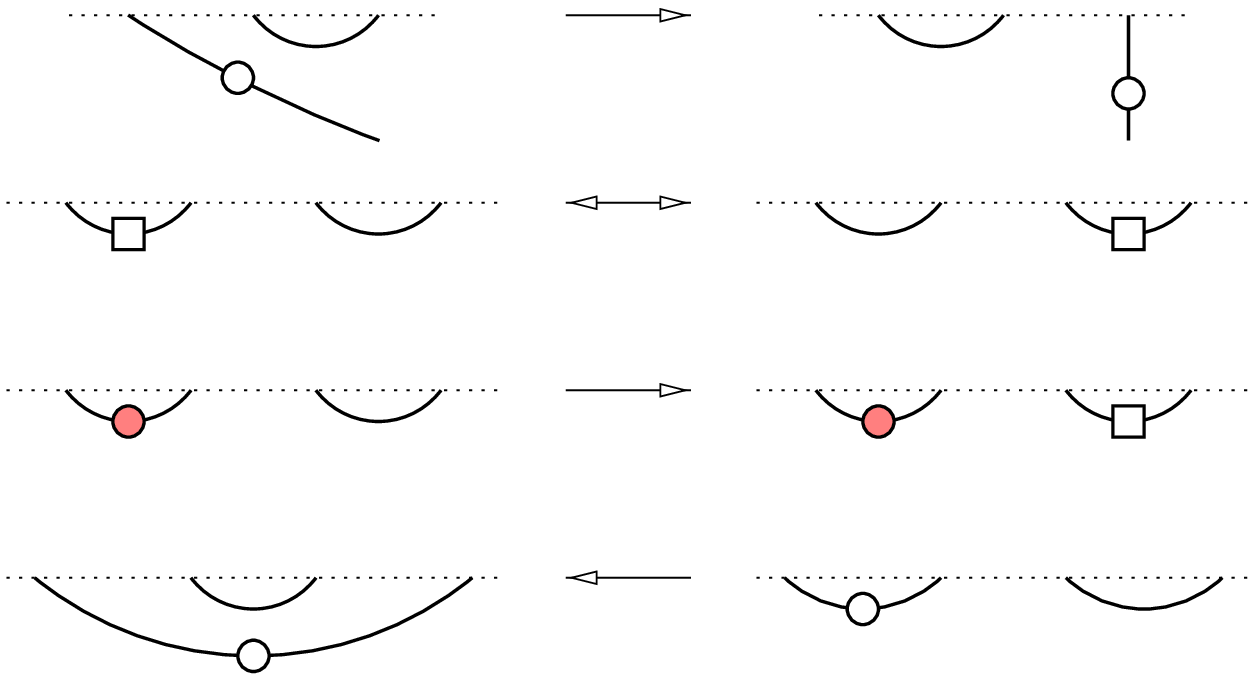}
}
\hfill}
}

The following property of the diagram calculus is analogous to that
described in Lemma 4.2.3 for type $H$.

\proclaim{Lemma 5.1.8}
Each element of the set $C_n$ of Theorem
2.2.5 may be obtained from the identity
element $1$ by repeated application of the following four procedures:
\item{\rm(i)}{left multiplication by $b_i$;}
\item{\rm(ii)}{right multiplication by $b_i$;}
\item{\rm(iii)}{left multiplication of an element $D'$ by $b_{s'}b_s -
1$ where $b_{s'}D' = (v+v^{-1})D'$ and $\{s, s'\} = \{s_1, s_2\}$;}
\item{\rm(iv)}{right multiplication of an element $D'$ by $b_sb_{s'} -
1$ where $D'b_{s'} = (v+v^{-1})D'$ and $\{s, s'\} = \{s_1, s_2\}$;}
\endproclaim

\demo{Proof}
The argument is similar to the argument in \cite{{\bf 5}, \S3.2} for type
$H$.  We first note that the set of Theorem 2.2.5 consists, in type
$B_2$, of the elements $$
1, b_1, b_2, b_1b_2, b_2b_1, (b_1b_2 - 1)b_1, (b_2b_1 - 1)b_2
,$$ which agrees with the hypotheses.

The general case uses lemmas 5.1.4--5.1.7.  The results about $\iota$
mean that we may copy the argument of \cite{{\bf 5}, Proposition 3.2.8},
from which the conclusion follows.
\qed\enddemo

\proclaim{Proposition 5.1.9}
The $\A$-spans of the following sets are equivalent:
\item{\rm(i)}{$\{t_w : w \in W_c\}$,}
\item{\rm(ii)}{$\{b_w : w \in W_c\}$,}
\item{\rm(iii)}{the set $C_n$ in the statement of Theorem 2.2.5.}
\endproclaim

\demo{Proof}
The equivalence of the first two follows from the fact that they are
both $\A$-bases for $TL(X)$; see the remarks in Definition 1.2.4.
We prove that the sets in (ii) or (iii) have the
same $\A$-span, implicitly using the fact that the set in (i) spans an
$\A$-algebra.

First consider $b_w$ for $w \in W_c$.  Since $b_i \in C_n$,
repeated applications of Lemma 5.1.1
show that $b_w$ lies in the $\A$-span of the set in (iii).
Conversely, Lemma 5.1.8 shows that each element of the set in (iii) is
a polynomial (over $\zed$) 
in the elements $b_i$, which establishes the reverse inclusion.
\qed\enddemo

\subhead 5.2 Main results in type $B$ \endsubhead

Since the set $C_n$ of Theorem 2.2.5 produces the correct integral
form, the arguments of \S3 and \S4 may be easily adapted to type $B$.
In general, these arguments are subsets of the arguments in type $H$:
the difference is that strings $s_1s_2s_1s_2$ and $s_2s_1s_2s_1$ are
not allowed in elements of $W_c$, which eliminates several of the most
complicated cases.  For reasons of space, we describe only the
necessary changes to the proof, together with statements of the main
intermediate results.

Lemmas 3.1.1--3.1.3 have direct analogues for type $B$, replacing the
diagram basis in type $H$ with the set $C_n$.  The set $C_n$ leads to
the definition of a projection $\pi_B$ analogous to $\pi_H$.  
The definitions of internal and lateral are the same as in type $H$.
Cases (ii) and (iii) in Proposition 3.1.9 cannot actually occur in
practice, but this is neither immediate nor necessary to get the
argument to work.

The definition of the $f$-basis in type $B$ is formally identical to type
$H$ in \S3.2, so we end up with a set which is formally the same as a subset of
the $f$-basis in type $H$.

The results of \S3.3 adapt readily to type $B$, the main difference
being that the case dealt with by Lemma 3.3.6 (ii) is not needed.

Copying the results of \S3.4, we obtain the following analogue of
Theorem 3.4.3.

\proclaim{Theorem 5.2.1}
The basis $\{f_w : w \in W_c\}$ is the canonical basis of
$TL(B_{n-1})$ in the sense of Theorem 1.2.5.
\endproclaim

The results of \S4 also adapt easily; the main changes to make are the
omissions of some of the cases.

For Lemma 4.1.1, we replace (i) (b) with its counterpart for type $B$
given in \cite{{\bf 7}, Lemma 2.1.2 (ii)}.  The case in Lemma 4.1.1
(ii) (b) cannot occur and can be ignored; the same goes for the situation
in Corollary 4.1.2 (iii).  Similar remarks hold for Lemma 4.1.6 (ii)
and Lemma 4.1.8 (ii).  This allows us to prove an analogue of
Proposition 4.1.9 which verifies \cite{{\bf 7},
Conjecture 1.2.4} for Coxeter systems of type $B$.

\proclaim{Proposition 5.2.2}
The structure constants for the canonical basis for $TL(B_{n-1})$ 
lie in $\A^{\geq 0}$.
\endproclaim

The cases to be ignored in \S4.2 are Lemma 4.2.2 (ii) and Lemma 4.2.3 (v)
and (vi).

Theorem 2.2.5 follows. \qed

\subhead 5.3 Concluding remarks \endsubhead

It is natural to wonder whether short proofs of Theorem 2.1.3 and
Theorem 2.2.5 exist.  It turns out (although we will not pursue this
here) that these results are closely related to the following
hypothesis.

\proclaim{Hypothesis 5.3.1}
Consider a generalized Temperley--Lieb algebra, $TL(X)$ with $t$-basis
$\{\te_w : w \in W_c\}$.  Then there is a symmetric, anti-associative,
nondegenerate $\A$-bilinear form $\lan, \ran$ on $TL(X)$ with respect to
which $$
\lan \te_w, \te_x \ran = \d_{w, x} \mod v^{-1}\A^-
,$$ where $\d$ is the Kronecker delta.
\endproclaim

By ``anti-associative'' above, we mean that $\lan \te_s \te_w, \te_x \ran =
\lan \te_w, \te_s \te_x \ran$.

Note that there is an analogous hypothesis for Hecke algebras $\H(X)$,
but that the existence of such a form is well-known and almost trivial (set
$\lan T_w, T_x \ran := \d_{w, x}q^{\ell(w)}$).  It is not hard to establish
Hypothesis 5.3.1 in type $H$ from Theorem 2.1.3 and in type $B$ from
Theorem 2.2.5, but proving it directly is curiously difficult.
Conversely, given Hypothesis 5.3.1, we may characterize the canonical
basis up to sign using the bilinear form as in \cite{{\bf 11}, \S14.2} and
then show, using a short argument, that
the relevant diagram basis satisfies this characterization and some
mild positivity assumptions.  We will return to the consideration of
Hypothesis 5.3.1 in a subsequent paper.

\head Acknowledgements \endhead

I thank J. Losonczy for suggesting many improvements to an
early draft of this paper.  I also thank the referee for some helpful
comments.

\leftheadtext{}
\rightheadtext{}

\end
\Refs\refstyle{A}\widestnumber\key{GL2}
\leftheadtext{References}
\rightheadtext{References}

\ref\key{{\bf 1}}
\by P. Etingof and M. Khovanov
\paper Representations of tensor categories and Dynkin diagrams
\jour Int. Math. Res. Not. 
\vol 5 \yr 1995 \pages 235--247
\endref

\ref\key{{\bf 2}}
\by C.K. Fan and R.M. Green
\paper On the affine Temperley--Lieb algebras
\jour Jour. L.M.S.
\vol 60 \yr 1999 \pages 366--380
\endref

\ref\key{{\bf 3}}
\by J.J. Graham
\book Modular representations of Hecke algebras and related algebras
\publ Ph.D. thesis
\publaddr University of Sydney
\yr 1995
\endref

\ref\key{{\bf 4}}
\by R.M. Green
\paper Generalized Temperley--Lieb algebras and decorated tangles
\jour J. Knot Th. Ram.
\vol 7 \yr 1998 \pages 155--171
\endref

\ref\key{{\bf 5}}
\bysame
\paper Cellular algebras arising from Hecke algebras of type $H_n$
\jour Math. Zeit.
\vol 229 \yr 1998 \pages 365--383
\endref

\ref\key{{\bf 6}}
\by R.M. Green and J. Losonczy
\paper Canonical bases for Hecke algebra quotients
\jour Math. Res. Lett.
\vol 6 \yr 1999 \pages 213--222
\endref

\ref\key{{\bf 7}}
\bysame 
\paper A projection property for Kazhdan--Lusztig bases
\jour Int. Math. Res. Not.
\vol 1 \yr 2000 \pages 23--34
\endref

\ref\key{{\bf 8}}
\by J.E. Humphreys
\book Reflection Groups and Coxeter Groups
\publ Cambridge University Press
\publaddr Cambridge
\yr 1990
\endref

\ref\key{{\bf 9}}
\by V.F.R. Jones
\paper Planar Algebras, I
\miscnote preprint, 1998
\endref

\ref\key{{\bf 10}}
\by D. Kazhdan and G. Lusztig
\paper Representations of Coxeter groups and Hecke algebras
\jour Invent. Math. 
\vol 53 \yr 1979 \pages 165--184
\endref

\ref\key{{\bf 11}}
\by G. Lusztig
\book Introduction to Quantum Groups
\publ Birkh\"auser \publaddr Basel \yr 1993
\endref

\ref\key{{\bf 12}}
\by P. Martin and H. Saleur
\paper The blob algebra and the periodic Temperley--Lieb algebra
\jour Lett. Math. Phys.
\vol 30 (3)
\yr 1994 
\pages 189--206
\endref

\ref\key{{\bf 13}}
\by R. Penrose
\paper Angular momentum: an approach to combinatorial space-time
\inbook Quantum Theory and Beyond (E. Bastin, Ed.)
\publ Cambridge University Press
\publaddr Cambridge
\yr 1971
\pages 151--180
\endref

\ref\key{{\bf 14}}
\by J.R. Stembridge 
\paper On the fully commutative elements of Coxeter groups 
\jour J. Algebraic Combin.
\vol 5 
\yr 1996 
\pages 353--385
\endref

\ref\key{{\bf 15}}
\bysame
\paper Some Combinatorial Aspects of Reduced Words in
Finite Coxeter Groups
\jour Trans. Amer. Math. Soc.
\vol 349 \yr 1997 \pages 1285--1332
\endref

\endRefs
